\documentstyle[12pt,aps,floats]{revtex}

\newcommand{\be}{\begin{equation}}
\newcommand{\ee}{\end{equation}}
\newtheorem{theorem}{Theorem}
\newtheorem{lemma}{Lemma}
\newtheorem{corollary}{Corollary}
\begin{document}
%
\title
{Removal of the resolvent-like energy dependence from
interactions and invariant subspaces of a total Hamiltonian}
\author{A.K.Motovilov}

\address{Laboratory of Theoretical Physics, JINR,
141980 Dubna, Moscow region, Russia\\
E--address: motovilv@thsun1.jinr.dubna.su}

\maketitle

\bigskip

\bigskip

\begin{abstract}
The spectral problem $(A + V(z))\psi=z\psi$ is considered where the main
Hamiltonian $A$ is a self-adjoint operator of sufficiently arbitrary
nature. The perturbation  $V(z)=-B(A'-z)^{-1}B^{*}$  depends
on the energy $z$ as resolvent of another self-adjoint operator $A'$.
The latter is usually interpreted as Hamiltonian describing an
internal structure of physical system. The operator $B$ is assumed to
have a finite Hilbert-Schmidt norm.  The conditions are formulated
when one can replace the perturbation $V(z)$ with an
energy-independent ``potential'' $W$ such that the Hamiltonian $H=A
+W$ has the same spectrum (more exactly a part of spectrum) and the
same eigenfunctions as the initial spectral problem.  The Hamiltonian
$H$ is constructed as a solution of the non-linear operator equation
$H=A+V(H)$. It is established that this equation is closely connected
with the problem of searching for invariant subspaces of the
Hamiltonian
$
{\bf H}=\left[ \begin{array}{lr}  A     &   B  \\
                                    B^{*} &   A'    \end{array}\right].
$
The orthogonality and expansion theorems are proved for eigenfunction
systems of the Hamiltonian $ H=A + W $. Scattering theory is
developed for this Hamiltonian in the case where the  operator $A$ has
continuous spectrum.
\end{abstract}
\bigskip

\bigskip

AMS subject classification (1980) 47B99, 81G05.
\pacs{02.30.Tb, 03.65.Nk, 12.40.--y}

\bigskip

\noindent Published in  
J. Math. Phys., 1995, vol. 36 (12), pp. 6647-6664. 
LANL E-print {\tt funct-an/9606002}.

\section{ INTRODUCTION}\label{Intro}
Perturbations, depending  on   the   spectral
parameter  (usually energy of system)  arise
in  a  lot   of   physical
problems (see papers~\cite{Dashen}--\cite{IMU} and Refs. therein).
In  particular,  such  are  the  interaction
potentials between clusters formed by quantum
particles~\cite{Dashen}--\cite{NarodetskyPalanga}.

The  perturbations  of this type appear
typically~\cite{Dashen}--\cite{Kerbikov}, \cite{Pavlov84}--\cite{IMU}
as a result  of
dividing the Hilbert space ${\cal H}$ of physical system in two subspaces,
${\cal H}={\cal H}_1\oplus{\cal H}_2$.
The first one, say ${\cal H}_1$, is interpreted as a space of
``external'' (for example, hadronic) degrees of freedom. The second one,
${\cal H}_2$, is associated with an ``internal'' 
(for example, quark) structure
of the system. The Hamiltonian ${{\bf H}}$ of the system looks as a matrix,
\be
\label{twochannel}
{{\bf H}}=\left[
\begin{array}{cc}
                     A_1              &          B_{12}          \\
                     B_{21}           &          A_2
\end{array}
\right]
\ee
with $A_1$, $A_2$, the channel Hamiltonians (self-adjoint
operators) and
$B_{12}$, $B_{21}=B_{12}^{*}$, the coupling operators.
Reducing the spectral problem
${{\bf H}}U=zU$ , $U=\{ u_1 ,u_2 \}$ to the
channel ${\alpha}$ only one gets the spectral problem
\be
\label{ini}
[A_{\alpha}  +V_{\alpha}  (z)]u_{\alpha}  =zu_{\alpha}  ,
\quad {\alpha} =1,2,
\ee
where the perturbation
\be
\label{epot}
V_{\alpha}  (z)=-B_{\alpha\beta} (A_{\beta} -z)^{-1}B_{\beta\alpha} ,
\quad {\beta} \neq {\alpha},
\ee
depends on the spectral parameter $z$ 
as the resolvent $(A_{\beta} - z)^{-1}$ of
the Hamiltonian $A_{\beta}$.
In more complicated cases $V_{\alpha}(z)$ 
can include also linear terms in respect
with  $z$. Other types of dependency of the potentials
$V_{\alpha}(z)$ on the spectral parameter z give, in a general way,
the spectral problems~(\ref{ini}) with a complex spectrum.

The present paper is a continuation of the author's
works~\cite{MotJMPh91}--\cite{SPbWorkshop} devoted to a study of the
possibility to ``remove'' the energy dependence from  perturbations
of the type~(\ref{epot}). Namely, in~\cite{MotJMPh91}--\cite{SPbWorkshop}
we construct such new potential $W_{\alpha}$ 
that spectrum of the Hamiltonian
$H_{\alpha}=A_{\alpha}+W_{\alpha}$ 
is a part of the spectrum of the problem~(\ref{ini}).
At the same time, the respective eigenvectors of $H_{\alpha}$ become
also those for~(\ref{ini}). 
Hamiltonians $H_{\alpha}$ are found as solutions of the
non-linear operator equations
\be
\label{basic}
H_{\alpha}=A_{\alpha} + V_{\alpha} (H_{\alpha})
\ee
first appeared in the paper~\cite{Braun}
by  M.A.Braun in connection with consideration
of the quasipotential equation. The operator-value function
$V_{\alpha}(Y)$ of the operator variable 
$Y$, $Y:\, {\cal H}_{\alpha}\rightarrow{\cal H}_{\alpha},$ is defined
by us in such a way (see Sec.~\ref{BasicEquation}) 
that eigenvectors $\psi$ of $Y$,
$Y\psi =z\psi$, become automatically those for $V_{\alpha}(Y)$ and
$V_{\alpha}(Y)\psi=V_{\alpha}(z)\psi.$

Note that the problem of removal of the dependence on a spectral
parameter from perturbation is interesting not only in itself.
In particular, there is a rather conceptual question
(see for instance Refs.~\onlinecite{McKellarMcKay},
\onlinecite{Schmid} and~\onlinecite{MotJMPh91})
concerning a use of the two--body energy--dependent potentials
in few--body nonrelativistic scattering problems. Since the energies
of  pair subsystems are not fixed in the N--body (N$\geq 3$) problem,
a direct embedding of such potentials into the few--body Hamiltonian
is impossible. Thus, the replacements of the type (\ref{epot})
energy--dependent potentials with the respective new
potentials $W_\alpha$ could be considered as a way to overcome
this difficulty.

In Ref.~\onlinecite{MotJMPh91}, the problem of the removal of the
energy dependence from the type (\ref{epot}) perturbations was
considered in details when one of the operators $A_{\alpha}$ is the
Schr\"{o}dinger operator in $L_2({\bf R}^n)$ and another one has a
discrete spectrum only. The reports~\cite{Kazimierz},
\cite{SPbWorkshop}  announce the results concerning the
equations~(\ref{basic}) and properties of their solutions $H_{\alpha}$ in a
rather more general situation where the Hamiltonian ${{\bf H}}$ may be
rewritten in terms of a two-channel variant of the Friedrichs model
investigated by O.A.Ladyzhenskaya and L.D.Faddeev in
Ref.~\onlinecite{LadyzhFaddeev} and L.D.Faddeev in
Ref.~\onlinecite{Faddeev64}.  In the paper~\cite{IMU} the
method~\cite{MotJMPh91}--\cite{SPbWorkshop} is applied for
construction of  an effective cluster Hamiltonian for atoms adsorbed
by the metal surface.

In the present paper, we specify the assertions
from~\cite{Kazimierz}, \cite{SPbWorkshop} and give proofs for them.
Also, we pay attention to an important circumstance disclosing a
nature of solutions of the basic equations (\ref{basic}). Thing is
that the potentials $W_{\alpha}=V_{\alpha}(H_{\alpha})$ 
may be presented in the form
$W_{\alpha}=B_{\alpha\beta} Q_{\beta\alpha}$ 
where the operators $Q_{\beta\alpha}$ satisfy the stationary
Riccati equations~(\ref{QbasicSym}) (see Sec.~\ref{BasicEquation}).
Exactly the same equations arise in the method of construction of
invariant subspaces developed by V.A.Malyshev and
R.A.Minlos~\cite{MalyshevMinlos1},~\cite{MalyshevMinlos2} for a class
of self-adjoint operators in statistical physics.  It follows from
the results of~\cite{MalyshevMinlos1}, \cite{MalyshevMinlos2} that
operators $H_{\alpha}$, ${\alpha}=1,2,$ determine in fact,  parts of the
two-channel Hamiltonian ${{\bf H}}$ acting in corresponding invariant
subspaces (see Theorem~\ref{ThInvariant} and comments to it).

Recently, the author came to know about the
work~\cite{LangerAdamyan}
by  V.M.Adamjan and H.Langer studying
the operator-value functions written in our notation as
$F_{\alpha}(z)=z-A_{\alpha}\pm 
B_{\alpha\beta}(A_{\beta}-z)^{-1}B_{\beta\alpha}$. 
In particular Adamjan and Langer show in this work
that a subset of eigenvectors of $F_{\alpha}$ can be chosen
to form a Riesz basis in
${\cal H}_{\alpha}$.  There is a certain intersection 
of their results and ours from
Refs.~\onlinecite{MotJMPh91}--\onlinecite{SPbWorkshop}.
However the methods are different.

The paper is organized as follows.

In Sec.~\ref{Initial} we describe the Hamiltonian ${{\bf H}}$ as a
two-channel variant of the Friedrichs
model~\cite{LadyzhFaddeev}, \cite{Faddeev64}.
We suppose that both operators $A_{\alpha}$, ${\alpha}=1,2,$ 
may have continuous
spectrum. When properties of objects connected with this spectrum
(wave operators and scattering matrices) are considered in
following sections, the coupling operators 
$B_{\alpha\beta}$ in~(\ref{twochannel})
are assumed to be integral ones with kernels 
$B_{\alpha\beta}(\lambda,\mu)$, the
 H\"{o}lder functions in both variables $\lambda,\mu$.

In Sec.~\ref{BasicEquation} the equations~(\ref{basic}) are studied. As in
Refs.~\onlinecite{MalyshevMinlos1},~\onlinecite{MalyshevMinlos2}
we suppose that
spectra $\sigma(A_1)$ and $\sigma(A_2)$ of the operators $A_1$ and
$A_2$ are separated, ${\rm dist}\{ \sigma(A_1), \sigma(A_2) \}>0$.
Existence of solutions of Eqs.~(\ref{basic}) is established
only in the case where the Hilbert-Schmidt norms $\| B_{\alpha\beta}\|_2$
of the coupling operators $B_{\alpha\beta}$ satisfy 
the condition
$\| B_{\alpha\beta}\|_2 < \frac{1}{2}{\rm dist}
\{ \sigma(A_1),\sigma(A_2) \} $.

In Sec.~\ref{Eigenfunctions} the eigenfunctions systems
of the operators $H_{\alpha}$ are studied and  theorems of their
orthogonality and completeness are proved. We show here
in particular that spectrum of the Hamiltonian ${{\bf H}}$ is distributed
between the solutions $H_1=A_1 + B_{12} Q_{21}$ and
$H_2=A_2 + B_{21} Q_{12}$, $Q_{21}=-Q_{12}^{*}, $
of the basic equations~(\ref{basic})
in such a way that $H_1$ and $H_2$ have not ``common'' eigenfunctions
$U=\{u_1, u_2\} $ of ${{\bf H}}$: simultaneously, the component $u_1$ can
not be eigenfunction for $H_1$, and the component $u_2$,  for $H_2$.

In Sec.~\ref{InnerProduct} we introduce new inner products
in the Hilbert spaces ${\cal H}_{\alpha}$, ${\alpha}=1,2,$
making the Hamiltonians $H_{\alpha}$ self-adjoint.

In Sec.~\ref{ScatteringProblem} we give a non-stationary formulation
of the scattering problem for a system described by the Hamiltonian
$H_{\alpha}$ constructed in the case where $A_\alpha$ has continuous
spectrum.  We show that this formulation is correct and scattering
operator is exactly the same as in initial spectral problem.

\section{ INITIAL SPECTRAL PROBLEM AND TWO-CHANNEL HAMILTONIAN}
\label{Initial}
Let $A_{1}$ and $A_{2}$ be self-adjoint  operators  acting,
respectively, in  ``external'', ${{\cal H}} _{1}$,  and  ``internal'',
${{\cal H}} _{2}$, Hilbert spaces. 
We study the spectral problem~(\ref{ini})
with  perturbation $V_{\alpha}(z)$  given by~(\ref{epot}).
We suppose that 
$B_{\alpha\beta} \in {{\bf B}} ({{\cal H}} _{\alpha}  ,
{{\cal H}} _{\beta} )$ where
${{\bf B}} ({{\cal H}} _{\alpha}  ,{{\cal H}} _{\beta} )$ is the
Banach  space  of   bounded  linear
operators acting from ${{\cal H}} _{\alpha}  $ to ${{\cal H}} _{\beta} $.

Note that the method developed in the present paper works also in the
case of more general perturbations$^{\onlinecite{RfunfNaboko}}$
 $V_{\alpha}(z)=-{\cal R}_{\alpha}(z)$ containing linear terms,
\be
\label{Rfunction}
  {\cal R}_{\alpha}(z)= N_{\alpha} z + 
B_{\alpha\beta}(A_{\beta}-N_{\beta} z -z)^{-1}B_{\beta\alpha}
\ee
with  $N_{\alpha}$, self-adjoint bounded operator
in ${\cal H}_{\alpha}$ such that 
$N_{\alpha}\geq ({\delta} -1)I_{\alpha}$ where
$\,\,{\delta}  > 0 $ and $I_{\alpha}$ 
is the identity operator in ${\cal H}_{\alpha}$.
Thing is that the equation~(\ref{ini}) with
$V_{\alpha}(z)=-{\cal R}_{\alpha}(z)$ 
can be easily rewritten in the form
(\ref{ini}),(\ref{epot}). To do this, 
one has only to make the replacements
$u_{\alpha}\rightarrow u_{\alpha}'=
(I_{\alpha}+N_{\alpha})^{1/2}u_{\alpha},$
$
\,\, A_{\alpha}\rightarrow A_{\alpha}'=
(I_{\alpha}+N_{\alpha})^{-1/2} A_{\alpha} (I_{\alpha}+N_{\alpha})^{-1/2}
$
and
$
  B_{\alpha\beta}\rightarrow B_{\alpha\beta}'=
(I_{\alpha}+N_{\alpha})^{-1/2}B_{\alpha\beta} 
(I_{\beta}+N_{\beta})^{-1/2}
$.
Therefore we shall consider further only the initial spectral problem
(\ref{ini}),(\ref{epot}).

We shall assume that operators $A_{\alpha}$, ${\alpha}=1,2,$ 
may have continuous
spectra $\sigma_{\alpha}^c$. To deal with 
these spectra we accept below some
presuppositions in respect with $A_{\alpha}$ 
restricting us to the case of
a two-channel variant of the Friedrichs
model~\cite{LadyzhFaddeev}, \cite{Faddeev64}.
Note that these presuppositions are not necessary
for a part of statements (Lemma~\ref{LQdefinition},
Theorems~\ref{ThSolvability}--\ref{ThDiscSpectrum}
and~\ref{ThInnerProduct}) which stay
correct also in general case.

The presuppositions are following.

At first, we assume that
Hamiltonian $H$ is defined  in  that  representation  where
operators $A_{\alpha}  , {\alpha} =1,2,$ are diagonal. We suppose
that continuous  spectra $\sigma ^{c}_{\alpha}  $
of  the  operators $A_{\alpha}  ,
{\alpha} =1,2,$ are absolutely continuous and consist of a finite
number of finite (and may be one or two infinite) intervals
$(a_{\alpha}^{(j)} ,b_{\alpha}^{(j)} )$,
$ -\infty\leq a_{\alpha}^{(j)}  < b_{\alpha}^{(j)}  \leq +\infty$,
$j=1,2,\ldots,n_{\alpha},$ $n_{\alpha} < \infty$.
At second, we suppose that discrete spectra $\sigma_{\alpha}^d$ of
the operators $A_{\alpha}$, ${\alpha}=1,2,$ 
do not intersect with $\sigma_{\alpha}^c$,
$\sigma_{\alpha}^d \bigcap\sigma_{\alpha}^c =\emptyset$,
and consist of a finite number of points with finite multiplicity.
In  this  case  the
space ${{\cal H}} _{\alpha}  $  may be present as the direct
integral~\cite{BirmanSolomiak}
\be
\label{Neumann}
{{\cal H}} _{\alpha}
=\sum\mbox{\hskip-2em}\int\limits_{\lambda\in\sigma_{\alpha} }
\oplus {{\cal G}} _{\alpha}  (\lambda )d\lambda
\equiv \mathop{\oplus}\limits_{\lambda\in\sigma^{d}_{\alpha}  }
{{\cal G}} _{\alpha}  (\lambda )\oplus
\int\limits_{\lambda\in\sigma_{\alpha}^c}
\oplus{{\cal G}} _{\alpha}  (\lambda )d\lambda , \quad
\sigma_{\alpha} =\sigma^{c}_{\alpha}\bigcup
\sigma^{d}_{\alpha}\subset{\bf R} .
\ee
The space ${{\cal H}} _{\alpha}  $
consists of the
measurable functions $f_{\alpha}  $
which are defined  on $\sigma _{\alpha}  $
and
have the values
$f_{\alpha}  (\lambda )$
from  corresponding  Hilbert
spaces               ${{\cal G}} _{\alpha}  (\lambda )$.
By $\langle\,\cdot\, ,\,\cdot\,\rangle$
we denote the inner  product  in
${{\cal H}} _{\alpha}  ,\qquad $
$$
\langle f_{\alpha}  ,g_{\alpha}  \rangle=
\mathop{\mbox{$\sum$\hskip-1em{\LARGE $\int$}\hskip0.1em}}
\limits_{\lambda\in\sigma_{\alpha}  }
(f_{\alpha}  (\lambda ),g_{\alpha}  (\lambda ))
\equiv \sum\limits_{\lambda\in\sigma^{d}_{\alpha}  }
(f_{\alpha}  (\lambda ),g_{\alpha}  (\lambda ))
+\int\limits_{\lambda\in\sigma_{\alpha}^c}
d\lambda (f_{\alpha}  (\lambda ),g_{\alpha}  (\lambda )),
$$
where $(\,\cdot\, ,\,\cdot\,)$ stands for inner product in
${{\cal G}} _{\alpha}  (\lambda )$. By  $|\cdot| $   we  denote
  norm  of  vectors  and
operators in ${{\cal G}} _{\alpha}  (\lambda )$ and by
$\parallel\cdot\parallel $,
the  norm  in ${{\cal H}} _{\alpha}  $.
 Operator $A_{\alpha}  $
acts in ${{\cal H}} _{\alpha}  $
as  the  independent  variable multiplication operator,
\be
\label{multi}
(A_{\alpha}  f_{\alpha}  )(\lambda )=
\lambda \cdot f_{\alpha}  (\lambda ), \quad {\alpha} =1,2.
\ee
It's domain
${\cal D}(A_{\alpha}  )$
consists of  those  functions $f_{\alpha}  \in {{\cal H}} _{\alpha}  $
which satisfy the condition
\newline
$\mathop{\mbox{$\sum$\hskip-1em{\LARGE $\int$}\hskip0.1em}}
\limits_{\lambda\in\sigma _{\alpha}  }\lambda ^{2}$
$|f_{\alpha}  (\lambda )|^{2}<\infty $. For  the
sake of simplicity we assume that
${{\cal G}} _{\alpha}(\lambda )$ does not depend on
$\lambda\in\sigma_{\alpha}^c$, i.e.
$
{{\cal G}} _{\alpha}(\lambda )\equiv {{\cal G}} ^{c}_{\alpha}
$
for each
$
\lambda \in \sigma^{c}_{\alpha}  .
$
Hence,
$
\int\limits_{\sigma ^{c}_{\alpha}  }
\oplus {{\cal G}} _{\alpha}  (\lambda )d\lambda =
L_{2}(\sigma ^{c}_{\alpha}  ,{{\cal G}} ^{c}_{\alpha}  )
\equiv {{\cal H}} ^{c}_{\alpha}
$.
By $E_{\alpha}  (d\lambda )$  we  denote  a
spectral   measure~\cite{BirmanSolomiak}   of
the operator $A_{\alpha}  ,$
$
A_{\alpha}  =\int\limits_{\sigma_{\alpha}}\lambda 
E_{\alpha}  (d\lambda )
$.
In the diagonal representation  considered,
the spectral projector $E_{\alpha}  $ acts on 
$f\in {{\cal H}} _{\alpha}  $ as
\be
\label{mera}
(E_{\alpha}  (\Delta  )f)(\lambda )=\chi_{\Delta  }(\lambda )f(\lambda )
\ee
for any Borelian set
$
\Delta  \subset \sigma _{\alpha}
$.
Here, $\chi _{\Delta  }$
is a characteristic function of
$\Delta  ,$
$ \chi _{\Delta  }(\lambda )=1$
if
$\lambda \in \Delta  $,
and $\chi _{\Delta  }(\lambda )=0$
if
$\lambda \not\in \Delta  $.

Let ${{\cal B}} ^{{\alpha}{\beta}}_{\theta\gamma}$
be a class of functions $F$  defined  on
$\sigma _{\alpha}  \times \sigma _{\beta} ,$ ${\alpha},{\beta}=1,2$,
for    each $\lambda \in \sigma _{\alpha}  , \mu \in \sigma _{\beta} $
as   operator
$F(\lambda ,\mu ):$
${{\cal G}} _{\beta} (\mu )\rightarrow {{\cal G}} _{\alpha}  (\lambda )$,
with $\parallel F\parallel_{{{\cal B}} } < \infty $,  where
$$
\parallel F\parallel _{{{\cal B}} } =
\mathop{\rm sup}\limits_{
\mbox{
\scriptsize
$\begin{array}{c}
\mu\in\sigma_{\beta}\\
\lambda\in\sigma_{\alpha}
\end{array}$
}
}(1+|\lambda |)^{\theta }
(1+|\mu |)^{\theta }| F(\lambda ,\mu ) |  +
$$
$$+\mathop{\rm sup}\limits_{
\mbox{
\scriptsize
$\begin{array}{c}
\lambda,\lambda'\in\sigma_{\alpha}^c\\
\mu\in\sigma_{\beta}
\end{array}$
}
}
\left\{
(1+|\mu |)^{\theta}
\frac{|F(\lambda ,\mu )-F(\lambda ',\mu )|}
{|\lambda  -\lambda '|^{\gamma }}
\right\}+
$$
$$
+\mathop{\rm sup}\limits_{
\mbox{
\scriptsize
$\begin{array}{c}
\lambda\in\sigma_{\alpha}\\
\mu,\mu'\in\sigma_{\beta}^c
\end{array}$
}
}
\left\{   (1+|\lambda|)^{\theta }
\frac{|F(\lambda ,\mu )-F(\lambda ,\mu ')|}{|\mu-\mu'|^{\gamma }}
\right\}+
$$
$$
+\mathop{\rm sup}\limits_{
\mbox{
\scriptsize
$\begin{array}{c}
\lambda,\lambda'\in\sigma_{\alpha}^c\\
\mu,\mu'\in\sigma_{\beta}^c
\end{array}$
} }
\left\{ \frac{|F(\lambda,\mu)-F(\lambda',\mu)-
F(\lambda,\mu')+F(\lambda',\mu')|}
{|\lambda-\lambda'|^{\gamma}|\mu-\mu'|^{\gamma }}
\right\} .
$$
With the norm $\parallel\cdot\parallel_{{{\cal B}} }$
this class will constitute a Banach space.
We introduce also the Banach space
${\cal M}_{\theta \gamma }(\sigma _{\alpha}  )$  of
functions $f$ defined on $\sigma _{\alpha}  $ with the norm
$$
\parallel f\parallel_{{\cal M}}
=\mathop{\rm sup}\limits_{\lambda\in\sigma_{\alpha}}
(1+|\lambda |)^{\theta }|f(\lambda)|
+\mathop{\rm sup}\limits_{\lambda,\lambda'\in\sigma_{\alpha}^c}
\frac{|f(\lambda ) - f(\lambda ')|}
{|\lambda  -\lambda '|^{\gamma }}<\infty .
$$
 The  value $f(\lambda )$
of  the  function $f\in {\cal M}_{\theta \gamma }(\sigma _{\alpha}  )$
is  an operator in ${{\cal G}} _{\alpha}  (\lambda )$.

Let $B_{\alpha\beta} $ be an integral  operator  with a
kernel $B_{\alpha\beta} (\lambda ,\mu )$
from the space ${{\cal B}} ^{{\alpha} {\beta} }_{\theta \gamma },$
$\theta >{1\over 2}, {1\over 2}<\gamma <1.$  We assume that
$B_{\alpha\beta} (\lambda ,\mu )$ is
 a compact operator, $B_{\alpha\beta} (\lambda ,\mu ):$
${{\cal G}} _{\beta} (\mu )\rightarrow {{\cal G}} _{\alpha}  (\lambda )$,
for each $\lambda \in \sigma _{\alpha}  , \mu \in \sigma _{\beta} $ and
 $B_{\alpha\beta} (\lambda ,\mu )=0$
if $\lambda $ belongs to the boundary of $\sigma ^{c}_{\alpha}  $
or $\mu $ belongs to the
boundary of $\sigma ^{c}_{\beta} $.

With this presuppositions the Hamiltonian ${{\bf H}}$
may be considered as a two-channel variant of the Friedrichs
model~\cite{LadyzhFaddeev}, \cite{Faddeev64}. Investigation of ${{\bf H}}$
repeats almost literally the analysis from
Ref.~\onlinecite{Faddeev64}. Therefore we describe
here only final results which are quite analogous to
\cite{LadyzhFaddeev}, \cite{Faddeev64}. These results are
following.

The  operator ${{\bf H}}$  is  self-adjoint  on  the  set
${\cal D}({{\bf H}})={\cal D}(A_{1})\oplus {\cal D}(A_{2})$.
Continuous  spectrum  of ${{\bf H}}$  is
situated on the set $\sigma _{c}({{\bf H}})=
\sigma ^{c}_{1}\cup \sigma ^{c}_{2}$.
Let ${{\bf H}}^{c}$ be the part of ${{\bf H}}$
acting in the invariant  subspace  corresponding  to
continuous  spectrum.   The  operator ${{\bf H}}^{c}$
is  unitary equivalent  to   the   operator
${{\bf H}}_{0}=A^{(0)}_{1}\oplus A^{(0)}_{2}$   with
$A^{(0)}_{\alpha},$ ${\alpha}=1,2,$ the restriction of the operator
$A_{\alpha}$ on ${\cal H}^{c}_{\alpha}$.
Namely,  there   exist  wave  operators $U^{(+)}$   and $U^{(-)}, $
$
U^{(\pm )}=
\left(
\begin{array}{lr}
u^{(\pm )}_{11}    &     u^{(\pm )}_{12}  \\
u^{(\pm )}_{21}    &     u^{(\pm )}_{22}
\end{array}
\right)=
$
$
s\!-\!\mathop{\rm lim}\limits_{t\rightarrow\mp\infty}
{\rm e}^{i{{\bf H}}t}{\rm e}^{-i{{\bf H}}_0 t},
$
 with  the   following   properties:
${{\bf H}}U^{(\pm )}=U^{(\pm )}{{\bf H}}_{0},$ 
$U^{(\pm )*}U^{(\pm )}=I,$
$U^{(\pm )}U^{(\pm )*}=  I-P.$
Here, $P$ is an
orthogonal projector on  subspace
corresponding to the discrete spectrum $\sigma _{d}({{\bf H}})$
of the operator ${{\bf H}}$.

The kernel $u^{(\pm )}_{{\alpha} {\alpha} }(\lambda ,\lambda ')$
of the operator $u^{(\pm )}_{{\alpha} {\alpha} },{\alpha} =1,2,$
represents an eigenfunction of the continuous spectrum
of the problem~(\ref{ini}) for
$
z=\lambda '\pm i0,
$
$
\lambda '\in \sigma ^{c}_{\alpha}
$, and  satisfies the integral equation
\be
\label{wavefunction}
u^{(\pm )}_{{\alpha} {\alpha} }(\lambda ,\lambda ')=
I_{\alpha}^c  {\delta}  (\lambda -\lambda ')-
[(A_{\alpha}  -\lambda '\mp i0)^{-1}
V_{\alpha}  (\lambda '\pm i0)u^{(\pm )}_{{\alpha} {\alpha} }]
(\lambda ,\lambda '),
\ee
where $I_{\alpha}^c  $
is identity operator  in
$
{{\cal G}} ^{c}_{\alpha}  , \lambda \in \sigma _{\alpha}
$.
For each concrete  sign  (plus  or  minus)  and  for  each
$
\lambda '\in \sigma ^{c}_{\alpha}, 
\lambda '\not\in \sigma _{d}({{\bf H}})
$
the function
$
u^{(\pm )}_{{\alpha} {\alpha} }(\lambda ,\lambda ')
$
is  an  unique solution of eq.(\ref{wavefunction})
in the class of the distributions
\be
\label{amplitude}
f^{(\pm )}_{\alpha}  (\lambda )=
I_{\alpha}  {\delta}  (\lambda -\lambda ')+
{f(\lambda )\over \lambda -\lambda '\mp i0},
\quad f\in {\cal M}_{\theta '\gamma '},
\ee
where
$
{1\over 2}<\theta '<\theta ,
$
$
{1\over 2}<\gamma '<\gamma
$.
At the same time
$$
u^{(\pm )}_{\alpha\beta} (\lambda ,\lambda ')=
-[(A_{\alpha}  -\lambda '\mp i0)^{-1}
B_{\alpha\beta} u^{(\pm )}_{{\beta} {\beta} }](\lambda ,\lambda '),
\quad {\beta} \neq {\alpha} ,
$$
is the  problem ~(\ref{ini})  eigenfunction  corresponding  to
$\lambda '\in \sigma_{\beta}^c$.

The functions $u_{{\beta}{\alpha}}^{(\pm)}$, 
${\alpha},{\beta}=1,2,$ can be explicitly
expressed in terms of kernels of the operator
$$
T(z)=B-B({{\bf H}}-z)^{-1}B,\quad
B=\left[\begin{array}{lr}   0       &   B_{12}  \\
                           B_{21}   &     0     \end{array}\right].
$$
Corresponding formulae read as
$$
u_{{\beta}{\alpha}}^{(\pm)}(\lambda,\lambda')=
{\delta} _{{\beta}{\alpha}}I_{\alpha}^c{\delta} (\lambda-\lambda')-
\frac{T_{{\beta}{\alpha}}(\mu,\lambda',\lambda'\pm i0)}
{\mu-\lambda'\mp i0}, \quad
\mu\in\sigma_{\beta},\,\, \lambda'\in\sigma_{\alpha}^c,
$$
with $t$-matrices
$$
T_{{\alpha}{\alpha}}=B_{\alpha\beta}
\left[z-A_{\beta}+B_{{\beta}{\alpha}}
(A_{\alpha}-z)^{-1}B_{\alpha\beta}\right]^{-1}
B_{{\beta}{\alpha}}
$$
      and
$$
T_{{\beta}{\alpha}}=B_{{\beta}{\alpha}}
\left[z-A_{\alpha}+B_{{\alpha}{\beta}}(A_{\beta}-z)^{-1}
B_{{\beta}{\alpha}}\right]^{-1}(z-A_{\alpha})=
$$
$$
=(z-A_{\beta})\left[z-A_{\beta}+
B_{{\beta}{\alpha}}(A_{\alpha}-z)^{-1}B_{\alpha\beta}\right]^{-1}
B_{{\beta}{\alpha}},\quad {\beta}\neq{\alpha}.
$$
Considering the equation for $T(z)$,
$T(z)=B-B(A-z)^{-1}T(z)$, $A=A_1\oplus A_2$, one shows
in the same way as in~\cite{LadyzhFaddeev}, \cite{Faddeev64}
that for all $z\in{\bf C}\setminus\sigma({{\bf H}})$, each kernel
$ T_{{\beta}{\alpha}}(\mu,\lambda,z)$, ${\alpha},{\beta}=1,2,$
belongs to the class
${{\cal B}} _{\theta'\gamma'}^{{\beta}{\alpha}}$
with arbitrary $\theta',\gamma'$
such that $\frac{1}{2}<\theta'<\theta$,
$\frac{1}{2}<\gamma'<\gamma$.
In respect with variable $z$,
the kernel of $T_{{\beta}{\alpha}}(z)$
is continuous in the 
${{\cal B}} _{\theta'\gamma'}^{{\beta}{\alpha}}$--norm
right up to the upper and lower borders of the set
$\sigma_c({{\bf H}})\setminus\sigma_d({{\bf H}})$.

Scattering operator $S=U^{(-)*}U^{(+)}$ for a system described by the
Hamiltonian ${{\bf H}}$ is unitary 
in ${\cal H}_{\alpha}^c$. It's kernels
$s_{{\beta}{\alpha}}(\mu,\lambda)$, 
${\alpha},{\beta}=1,2$, are given by expressions
\be
\label{scattering}
  s_{{\beta}{\alpha}}(\mu,\lambda)={\delta} (\mu-\lambda)
  \left[       {\delta} _{{\beta}{\alpha}}I_{\alpha}^c-
  2\pi i\, T_{{\beta}{\alpha}}(\mu,\lambda,\lambda+i0)  \right].
\ee

By $U_{j},$ $j=1,2,\ldots$,
we  denote  eigenvectors,
$U_{j}=\{u^{(j)} _{1},u^{(j)} _{2}\},$
$ U_{j}\in {\cal D}({{\bf H}})$,
$\parallel U_{j} \parallel =1,$ and by 
$z_{j},$ $ z_{j}\in {\bf R}$, the
respective eigenvalues of  the   operator ${{\bf H}}$  discrete
spectrum $\sigma_d({{\bf H}})$. The component
$u^{(j)} _{\alpha}  ,$ ${\alpha} =1,2,$ of the  vector $U_{j}$
is a solution of  Eq.~(\ref{ini})  at $z=z_{j}$.
If $z_{j}\!\in\!\sigma ^{c}_{\beta} $  then
$(B_{\beta\alpha} u^{(j)} _{\alpha}  )(z_{j})=0.$
%
\section{ CONSTRUCTION OF THE OPERATORS $H_\alpha$}
\label{BasicEquation}
%
The paper is  devoted  to  construction  of  such
operator $H_{\alpha}  $
that it's each eigenfunction
$u_{\alpha}  , H_{\alpha}  u_{\alpha}  =zu_{\alpha}  $,
together with eigenvalue $z$,  satisfies  Eq.~(\ref{ini}).  This
operator will be found as a solution of the non-linear
operator equation~(\ref{basic}). To obtain this equation we need
the following operator-value function $V_{\alpha}  (Y)$
of the  operator
variable $Y:\quad$
$$
V_{\alpha}  (Y)=B_{\alpha\beta}\int\limits_{\sigma _{\beta} }
E_{\beta} (d\mu )B_{\beta\alpha} (Y-\mu )^{-1},
$$
$Y:$ ${{\cal H}} _{\alpha}  \rightarrow {{\cal H}} _{\alpha}  $ .
We suppose here that
$
(Y-\mu I)^{-1}\in L_{\infty }
(\sigma _{\beta} ,{{\bf B}}({{\cal H}}_{\alpha},{{\cal H}}_{\alpha}))
$
if $\mu \in \sigma _{\beta} $.
This means that $\sigma_{\beta}$ has not 
to be included into the spectrum of
the operator $Y$. The integral
$
Q(T)=\int\limits_{\sigma _{\beta} }
E_{\beta} (d\mu )B_{\beta\alpha} T(\mu )
$
may be constructed in the same way as integrals of scalar functions
over spectral measure~\cite{BirmanSolomiak}.
Namely, we consider $Q(T)$ for
$
T\in L_{\infty }(\sigma _{\beta} ,{{\bf B}} ({{\cal H}} _{\alpha}  ,
{{\cal H}} _{\alpha}  )),
$
$
\|{T}\|_{\infty }=
$
$
E_{\beta}\!-\!\mathop{\rm sup}\limits_{\mu\in\sigma_{\beta}}
\|T(\mu)\|<\infty
$
(hereafter $E_{\beta}\!-\!\mathop{\rm sup}\limits_{}$ 
means the supremum with respect to $E_\beta$),
as a limit value, with respect to the operator
norm in ${{\bf B}}({{\cal H}}_{\alpha},{{\cal H}}_{\alpha})$,
of finite integral sums for piecewise--constant 
operator--value  functions
approximating $T$  in
$
L_{\infty }(\sigma _{\beta} ,{{\bf B}}
({{\cal H}} _{\alpha}  ,{{\cal H}} _{\alpha}  ))
$.
We show the existence of this integral at least in the case where
the Hilbert-Schmidt  norm $\|{B}_{{\alpha}{\beta}}\|_{2}$ is finite.

\begin{lemma}\label{LQdefinition}
Let
$
T\in L_{\infty }(\sigma _{\beta} ,
{{\bf B}} ({{\cal H}} _{\alpha}  ,{{\cal H}} _{\alpha}  ))
$
 and
$
\parallel B_{\alpha\beta} \parallel_{2}<\infty $.
 Then  the
integral $Q(T)$   exists  as  a   bounded   operator in ${\cal H}_\alpha$
with norm satisfying the estimate
$
\|Q(T)\|\le \| T\|_{\infty }\!\cdot\!\|B_{\beta\alpha}\|_{2}
$.
\end{lemma}

\noindent{\sl Proof}. We prove the  Lemma  in  the
diagonal  representation~(\ref{Neumann}),~(\ref{multi}).
By~(\ref{mera}) we have
$$
(Qf)(\mu )=
\mathop{\mbox{$\sum$\hskip-1em{\LARGE $\int$}\hskip0.1em}}
\limits_{\sigma_{\alpha}} B_{\beta\alpha} (\mu ,\lambda )
(T(\mu )f)(\lambda) d\lambda
$$
for
any $f\in {{\cal H}} _{\alpha}  $.
It means that
$$
| (Qf)(\mu )|^{2}\!\le\!\sum
\mbox{\hskip-2em}\int\limits_{\lambda\in\sigma_{\alpha}}
d\lambda \cdot
| B_{\beta\alpha} (\mu,\lambda )|^{2}
\sum\mbox{\hskip-2em}\int\limits_{\lambda\in\sigma_{\alpha}}d\lambda 
|({T(\mu)}f)(\lambda)|^{2}=
$$
$$
=\sum\mbox{\hskip-2em}
\int\limits_{\lambda\in\sigma_{\alpha}} 
d\lambda|B_{{\beta}{\alpha}}(\mu,\lambda)|^{2}\cdot
\|T(\mu)f\|^2
\leq \sum\mbox{\hskip-2em}
\int\limits_{\lambda\in\sigma_{\alpha}}
d\lambda|B_{{\beta}{\alpha}}(\mu,\lambda)|^{2}\cdot
\|T(\mu)\|^{2}\cdot\|f\|^2.
$$
Hence,
integrating  over
$\mu \in \sigma_{\beta} $
we  come  to  the   relation
$$
\parallel{Qf}\parallel^{2}\le
\parallel{B_{{\beta}{\alpha}}}\parallel^2_2\!\cdot\!
\parallel{T}\parallel^{2}_{\infty }\!\cdot\!
\parallel{f}\parallel^{2}
$$
which completes the proof.
\medskip

Let us suppose that
$
(H_{\alpha}  -\mu I)^{-1}\in L_{\infty }
(\sigma _{\beta} ,{{\bf B}} 
({{\cal H}} _{\alpha}  ,{{\cal H}} _{\alpha}  ))
$.
We note that if
$H_{\alpha}  \psi_{\alpha}  =z\psi_{\alpha}  $,
then automatically
$$
V_{\alpha}  (H_{\alpha}  )\psi_{\alpha}  =
B_{\alpha\beta}\int\limits_{\sigma _{\beta} }
E_{\beta} (d\mu )
B_{\beta\alpha} (z-\mu )^{-1}\psi_{\alpha}  =
$$
\be
\label{Vaction}
=B_{\alpha\beta} (z-A_{\beta} )^{-1}B_{\beta\alpha} \psi_{\alpha}
=V_{\alpha}  (z) \psi_{\alpha}  .
\ee
It follows from~(\ref{Vaction})  
that $H_{\alpha}  $  satisfies the  relation
$
H_{\alpha}  \psi_{\alpha}
=(A_{\alpha}  +V_{\alpha}  (H_{\alpha}  )) \psi_{\alpha}
$
and we can  spread  this  relation
over all the linear combinations of $H_{\alpha}  $ eigenfunctions.
Supposing that the  eigenfunctions  system  of $H_{\alpha}  $  is
dense in ${\cal H}_{\alpha}$ we spread this
equation over ${\cal D}(A_{\alpha}  )$. As a  result
we come to the desired {\it basic equation}
(\ref{basic})
for $H_{\alpha}  $ (see also Refs.~\onlinecite{Braun},
\onlinecite{MotJMPh91}--\onlinecite{SPbWorkshop}).
Eq.~(\ref{basic}) means that the construction 
of the operator $H_{\alpha}  $
comes to the searching for the operator
\be
\label{Qbasic}
Q_{\beta\alpha}
=\int\limits_{\sigma_{\beta}} E_{\beta} (d\mu )
B_{\beta\alpha} (H_{\alpha}  -\mu )^{-1}.
\ee
Since
$
H_{\alpha}  =A_{\alpha}  +B_{\alpha\beta} Q_{\beta\alpha}
$,
we have
\be
\label{Qbasic1}
Q_{\beta\alpha} =
\int\limits_{\sigma _{\beta} }E_{\beta} (d\mu )
B_{\beta\alpha} (A_{\alpha}  +
B_{\alpha\beta} Q_{\beta\alpha} -\mu )^{-1},
\quad {\beta} \neq {\alpha}.
\ee
In this paper we restrict ourselves to the  study
of Eq.~(\ref{Qbasic1})
solvability  only  in  the  case where
spectra $\sigma_{1}$ and $\sigma_{2}$
are  separated,
\be
\label{dist}
         d_{0}={\rm dist}(\sigma _{1},\sigma _{2})>0.
\ee
Using the Lemma~\ref{LQdefinition} and 
the contracting  mapping  Theorem,
we prove the following:
\begin{theorem}\label{ThSolvability}
 Let $M_{\beta\alpha} ({\delta}  )$
 be  a  set  of   bounded operators
$X,$ $X:$
${{\cal H}} _{\alpha}  \rightarrow {{\cal H}} _{\beta} $,
 satisfying the inequality
$ \|{X}\| $ $ \le{\delta}     $ with ${\delta}  > 0$.
 If this ${\delta} $ and the norm
$\|{B_{\alpha\beta} }\|_{2}$
 satisfy the condition
$
\| B_{\alpha\beta}\|_{2}< d_{0}
\mathop{\rm min}\limits_{} 
\{  \frac{1}{1+{\delta} },
\frac{{\delta} }{1+{\delta} ^2}  \}$,
 then Eq.~(\ref{Qbasic1})
 is  uniquely
 solvable in
$M_{\beta\alpha} ({\delta}  )$.
\end{theorem}

\noindent{\sl Proof.} Let
\be
\label{funF}
F(X)=\int\limits_{\sigma_{\beta}}
E_{\beta}(d\mu)B_{{\beta}{\alpha}}
(A_{\alpha}+B_{\alpha\beta} X-\mu)^{-1}
\ee
with $X$, the operator from 
${{\bf B}} ({\cal H}_{\alpha},{\cal H}_{\beta})$.

Firstly, consider conditions where the function $F$ maps the set
$M_{\alpha}({\delta} )$ into itself. We suppose here that 
$B_{\alpha\beta}$ and $X$ are such
that
\be
\label{ne1}
     \|B_{\alpha\beta}\|_2 \|X\| 
    \leq {\delta}  \|B_{\alpha\beta}\|_2 < d_0
\ee
and consequently, $\|B_{\alpha\beta} X\|\leq d_0 . $
This means that spectrum of the operator 
$A_{\alpha} +B_{\alpha\beta} X$
does not intersect with the set $\sigma_{{\beta}}$. 
Hence, the resolvent
$(A_{\alpha} +B_{\alpha\beta} X-\mu)^{-1}$
exists and is bounded for any $\mu\in\sigma_{\beta}$.
Thus, by Lemma~\ref{LQdefinition} we have
$$
  \|F(X)\| \leq \|B_{\alpha\beta}\|_2\cdot 
E_{\beta}\!-\!\mathop{\rm sup}\limits_{\mu\in\sigma_{\beta}}
  \|(A_{\alpha} +B_{\alpha\beta} X-\mu)^{-1}\|.
$$
Due to identity
$$
(A_{\alpha} +B_{\alpha\beta} X-\mu)^{-1}=
\left( I+(A_{\alpha}-\mu)^{-1}B_{\alpha\beta} X\right)^{-1}
  (A_{\alpha}-\mu)^{-1}
$$
and inequality $ \|B_{\alpha\beta}\| \leq \|B_{\alpha\beta}\|_2 $
we make estimation
$$
\| (A_{\alpha} +B_{\alpha\beta} X-\mu)^{-1} \|  \leq
\frac{1}{ 1 - \| (A_{\alpha}-\mu)^{-1}\| \|B_{\alpha\beta}\|_2 \|X\| }
\|(A_{\alpha}-\mu)^{-1}\|  \leq
$$
\be
\label{estires}
\leq \frac{1}{ 1-\frac{1}{d_0} \|B_{\alpha\beta}\|_2 {\delta}  }\cdot
\frac{1}{d_0}=
\frac{1}{d_0 - \|B_{\alpha\beta}\|_2 {\delta} }.
\ee
Therefore, the set $M_{\alpha}({\delta} )$ 
will be mapped by $F$ into itself if
 $\|B_{\alpha\beta}\|_2$ and ${\delta} $ are such that
\be
\label{ne2}
 \|B_{\alpha\beta}\|_2 \cdot 
\frac{1}{d_0 - \|B_{\alpha\beta}\|_2{\delta} } \leq {\delta} .
\ee

Secondly, study conditions for the 
function $F$ to be a contracting mapping.
Now, we consider the difference
$$
F(X)-F(Y)=\int\limits_{\sigma_{\beta}} 
E_{\beta}(d\mu) B_{\beta\alpha}
\left[ (A_{\alpha}+B_{\alpha\beta} X -\mu)^{-1} 
- (A_{\alpha}+B_{\alpha\beta} Y -\mu)^{-1} \right]=
$$
$$
=\int\limits_{\sigma_{\beta}} E_{\beta}(d\mu) B_{\beta\alpha}
 (A_{\alpha}+B_{\alpha\beta} X -\mu)^{-1} 
B_{\alpha\beta} (Y-X) (A_{\alpha}+B_{\alpha\beta} Y -\mu)^{-1}.
$$
Again, by Lemma~\ref{LQdefinition}, we have
$$
       \|F(X)-F(Y)\|  \leq
$$
$$
\leq \|B_{\alpha\beta}\|_2^2\cdot
\mathop{\rm sup}\limits_{\mu\in\sigma_{\beta}}
        \| (A_{\alpha}+B_{\alpha\beta} X -\mu)^{-1}\| \cdot
   \mathop{\rm sup}\limits_{\mu\in\sigma_{\beta}} 
\| (A_{\alpha}+B_{\alpha\beta} Y -\mu)^{-1}\| \cdot \|(Y-X)\| .
$$
With~(\ref{estires}) we come to the estimate
$$
\|F(X)-F(Y)\| \leq \|B_{\alpha\beta}\|_2^2 \cdot
\frac{1}{(d_0 -\|B_{\alpha\beta}\|_2{\delta} )^2}\cdot \|Y-X\|.
$$
The function $F$ becomes a contracting mapping if
\be
\label{ne3}
 \frac{\|B_{\alpha\beta}\|_2^2}
{(d_0 -\|B_{\alpha\beta}\|_2{\delta} )^2} < 1.
\ee
Solving system of the inequalities~(\ref{ne1}),~(\ref{ne2}) and
(\ref{ne3}) we find
$$
\|B_{\alpha\beta}\|_2 < d_0 \mathop{\rm min}\limits_{} \left\{
\frac{{\delta} }{1+{\delta} ^2}, \frac{1}{1+{\delta} } \right\}
$$
and this completes the proof of Theorem~\ref{ThSolvability}.

\begin{corollary}\label{CBestimate}
Equation (\ref{Qbasic1})
is uniquely solvable in the unit
ball $M_{\beta\alpha}(1)\subset$ 
${{\bf B}} ({\cal H}_{\alpha},{\cal H}_{\beta})$
for any $B_{\alpha\beta}$ such that
\be
\label{ineqB}
   \|B_{\alpha\beta}\|_2 < \frac{1}{2} d_0.
\ee
\end{corollary}

To prove the inequality~(\ref{ineqB}), note that
$
\mathop{\rm max}\limits_{{\delta} \geq 0}\,
\mathop{\rm min}\limits_{} \left\{
\frac{{\delta} }{1+{\delta} ^2}, \frac{1}{1+{\delta} } \right\}
=\frac{1}{2}
$
(at ${\delta} =1$). Hence, if~(\ref{ineqB}) 
takes place then the function
(\ref{funF}) is a contracting mapping of 
the unit ball $M_{\beta\alpha}(1)$ into
itself.
\medskip

\noindent{\bf Remark.} In the  proofs of Lemma~\ref{LQdefinition}
and Theorem~\ref{ThSolvability}
we did not use the assumption
about finiteness of the numbers $n_{\alpha}$ 
of intervals included in continuous
spectra $\sigma_{\alpha}^c$ of the operators 
$A_{\alpha}$, ${\alpha}=1,2$. Really,
these assertions take place in the case of arbitrary spectrum
$\sigma_{\alpha} $.

Finiteness at least of one of the numbers $n_1$ and $n_2$ will be used
at the moment.
If $n_1$ and/or $n_2$ are finite and
\be
\label{BQd}
\|B_{\alpha\beta} Q_{\beta\alpha}\| < d_0=
{\rm dist}\{\sigma_1,\sigma_2\},
\quad {\alpha}=1,2,\, {\beta}\neq{\alpha},
\ee
we can state  that
\be
\label{resolvest}
\|( A_{\alpha} + B_{\alpha\beta} Q_{\beta\alpha} -
\mu )^{-1}\| \leq \frac{C_{\alpha\beta}}{1+|\mu|}, \,\,
{\alpha}=1,2, \quad 
\mbox{at any $\mu\in\sigma_{\beta}$, ${\beta}\neq{\alpha}$,}
\ee
with some $C_{\alpha\beta} > 0$, 
$C_{\alpha\beta}\sim 1/(d_0 - \|B_{\alpha\beta} Q_{\beta\alpha}\| )$.
Of course this
estimate is essential only in the case 
where $\sigma_{\beta}$ is unbounded.
It follows immediately from Eq.~(\ref{Qbasic1}) that
if $n_1$ and/or $n_2$ are finite then
$Q_{\beta\alpha} f_{\alpha}\in 
{\cal D}(H_{\beta})={\cal D}(A_{\beta})$
for any $f_{\alpha}\in {\cal H}_{\alpha}$.

In this case we can rewrite Eq.~(\ref{Qbasic1}) in  symmetric
form as a stationary Riccati equation
(see the book~\cite{Reid} and Refs. therein),
\be
\label{QbasicSym}
Q_{\beta\alpha} A_{\alpha}  -
A_{\beta} Q_{\beta\alpha} +
Q_{\beta\alpha} B_{\alpha\beta} Q_{\beta\alpha} =B_{\beta\alpha} .
\ee
To make  this,  it  is  sufficient  to  calculate  the
expression
$
Q_{\beta\alpha} H_{\alpha}  -A_{\beta} Q_{\beta\alpha}
$
for both parts of eq.(\ref{Qbasic1}) having in mind that we apply it
to $f_{\alpha}\in {\cal D}(H_{\alpha})$. Did, we have
$$
Q_{\beta\alpha} H_{\alpha}  -A_{\beta} Q_{\beta\alpha} =
Q_{\beta\alpha}(A_{\alpha}+
B_{\alpha\beta} Q_{\beta\alpha})-A_{\beta} Q_{\beta\alpha}=
Q_{\beta\alpha} A_{\alpha} -A_{\beta} Q_{\beta\alpha} +
Q_{\beta\alpha} B_{\alpha\beta} Q_{\beta\alpha}.
$$
On the other hand,
$$
Q_{\beta\alpha} H_{\alpha}  -A_{\beta} Q_{\beta\alpha} =
\int\limits_{\sigma_{\beta}}[E_{\beta}(d\mu)
B_{\beta\alpha}(H_{\alpha}-\mu)^{-1}H_{\alpha}-
\mu E_{\beta}(d\mu)B_{\beta\alpha}(H_{\alpha}-\mu)^{-1}]=
B_{\beta\alpha}.
$$
One finds immediately from  Eqs.~(\ref{QbasicSym}),
${\alpha}=1,2$, that
if $Q_{\beta\alpha}$ gives solution 
$H_{\alpha}=A_{\alpha}+B_{\alpha\beta} Q_{\beta\alpha}$ 
of  the problem
(\ref{basic}) in the channel ${\alpha}$  then
\be
\label{Qadj}
Q_{\alpha\beta}=-Q_{\beta\alpha}^{*}=
-\int\limits_{\sigma_{\alpha} }( H^{*}_{\alpha}
-\mu)^{-1}B_{\alpha\beta} E_{\beta}(d\mu )
\ee
gives analogous solution 
$H_{\beta}=A_{\beta}+B_{\beta\alpha} Q_{\alpha\beta}$ 
in the channel ${\beta}$.

\begin{theorem}\label{ThInvariant}
Let $Q_{\beta\alpha}$, 
$Q_{\beta\alpha}\in{{\bf B}}({\cal H}_{\alpha},{\cal H}_{\beta})$,
be a solution of Eq.~(\ref{QbasicSym}) satisfying together with
$Q_{\alpha\beta}=-Q_{\beta\alpha}^{*}$ the conditions
(\ref{resolvest}). Then the transform
${{\bf H}}'={\cal Q}^{-1}{{\bf H}}{\cal Q}$ with
$
{\cal Q}=\left[\begin{array}{lr}   I_1      &   Q_{12} \\
                                   Q_{21}   &     I_2
\end{array}\right]
$
reduces the operator ${{\bf H}}$ to the block--diagonal form,
${{\bf H}}'=\mathop{\rm diag}\{ H_1, H_2 \}$ 
where $H_{\alpha}=A_{\alpha} + B_{\alpha\beta} Q_{\beta\alpha},$
${\alpha},{\beta}=1,2,$ ${\beta}\neq{\alpha}.$ 
At the same time, the operators
$
{\cal O}_{\alpha}=\left[\begin{array}{lr}   I_{\alpha}      &     0  \\
                                      Q_{\beta\alpha}     &     I_{\beta}
\end{array}\right]
$
reduce the Hamiltonian ${{\bf H}}$,
$
{{\bf H}}=\left[\begin{array}{lr}   
                          A_{\alpha}      &     B_{\alpha\beta}  \\
                      B_{\beta\alpha}     &     A_{\beta}
\end{array}\right],
$
to triangular form,
$
{{\bf H}}^{(\alpha)}\equiv 
{\cal O}_{\alpha}^{-1} {{\bf H}}{\cal O}_{\alpha}=
\left[\begin{array}{lr}     H_{\alpha}      &     B_{\alpha\beta}  \\
                              0       &     H_{\beta}^{*}
\end{array}\right].
$
\end{theorem}

\noindent{\sl Proofs} of both statements are done by direct substituting
 ${\cal Q}$ and ${\cal O}_{\alpha}$ into the definitions 
of ${{\bf H}}'$ and ${{\bf H}}^{(\alpha)}$ and
using the equations~(\ref{QbasicSym}) .

We have to note only that operator ${\cal Q}$ 
is invertible since, according to
~(\ref{Qadj}),
\be
\label{X}
X_{\alpha}=I_{\alpha}-Q_{\alpha\beta} Q_{\beta\alpha} = 
I_{\alpha} +Q_{\alpha\beta} Q_{\alpha\beta}^{*} \geq I_{\alpha}, 
\quad {\alpha}=1,2,
\ee
and
\be
\label{Qrevers}
{\cal Q}^{-1}=\left[  \begin{array}{lr}   X_1^{-1}   &   0 \\
                                        0    &     X_2^{-1}
\end{array}\right]   \cdot
\left[  \begin{array}{lr}    I_1             &      -Q_{12} \\
                           -Q_{21}         &         I_2
\end{array}\right] .
\ee

\begin{corollary}\label{CReducing}
Subspaces
${\cal H}^{(\alpha)} = {\cal O}  _{\alpha}
({\cal H}_{\alpha}\oplus\{  {\bf 0} \} ) =$
$
     \{ f:\,\, f=\{ f_{\alpha}, f_{\beta} \}\in{\cal H},\,\,
     f_{\alpha}\in{\cal H}_{\alpha},\,\, 
f_{\beta}=Q_{\beta\alpha} f_{\alpha} \}
$
are orthogonal, $  {\cal H}^{(1)} \perp {\cal H}^{(2)}, $
and reducing for ${{\bf H}}$,
${{\bf H}}\left( {\cal D}({{\bf H}})\bigcap
{\cal H}^{(\alpha)} \right)\subseteq {\cal H}^{(\alpha)}. $
\end{corollary}

Really, if $f\in{\cal H}^{(\alpha)}$, $g\in{\cal H}^{(\beta)}$ and
$f=\{  f_{\alpha}, Q_{\beta\alpha} f_{\alpha} \},$ 
$g=\{  Q_{\alpha\beta} g_{\beta}, g_{\beta} \},$ then
$\langle f,g \rangle = $
$\langle f_{\alpha}, Q_{\alpha\beta} g_{\beta} \rangle + 
\langle  Q_{\beta\alpha} f_{\alpha}, g_{\beta} \rangle = 0  $
since $Q_{\beta\alpha} = -Q_{\alpha\beta}^{*}$.
The invariance of ${\cal H}^{(\alpha)}$, ${\alpha}=1,2,$ 
in respect with ${{\bf H}}$ follows
 from the equality $ {{\bf H}}{\cal Q} ={\cal Q} {{\bf H}}'.$

Assertions quite analogous to the Theorem~\ref{ThInvariant}
and Corollary~\ref{CReducing} one
can find in Refs.~\onlinecite{MalyshevMinlos1},
\onlinecite{MalyshevMinlos2}.
Solvability (for sufficiently small $\|B_{\alpha\beta}\|$)
of the equation~(\ref{QbasicSym}) was proved
in~\onlinecite{MalyshevMinlos1},
\onlinecite{MalyshevMinlos2} by rather different method
also in the supposition~(\ref{dist}).
\medskip

\noindent{\bf Remark}.  It follows from
Theorem~\ref{ThInvariant}  that operator
$  \tilde{{\cal Q}}={\cal Q} X^{-1/2} $ 
with $X=\mathop{\rm diag}\{X_1,X_2\}$
is unitary. Consequently, the operator
 $ {{\bf H}}''= \tilde{{\cal Q}}^{*} {{\bf H}} 
\tilde{{\cal Q}}=$ $X^{1/2} {{\bf H}}' X^{-1/2}$
becomes self-adjoint in ${\cal H} $.
Since $ {{\bf H}}''=\mathop{\rm diag}\{  H''_1, H''_1  \}$
with $ H''_{\alpha}=X_{\alpha}^{1/2} H_{\alpha} X_{\alpha}^{-1/2}$,
 the operators $H''_{\alpha}$, $\,\,\,{\alpha}=1,2,$
are self-adjoint on ${\cal D}(A_{\alpha})$
in ${\cal H}_{\alpha} $. Moreover the operators
$
{{\bf H}}^{(\alpha)}=\tilde{{\cal Q}}\cdot 
\mathop{\rm diag}\{ H''_{\alpha}, 0\}\cdot  \tilde{{\cal Q}}^{*} =
$
${\cal Q} \cdot\mathop{\rm diag}\{ H_{\alpha}, 0 \}\cdot 
{\cal Q}^{-1} $  represent
parts of the Hamiltonian ${{\bf H}}$ in the 
corresponding invariant subspaces
${\cal H}^{(1)}$ and ${\cal H}^{(2)}$
(see also Refs.~\onlinecite{MalyshevMinlos1},
\onlinecite{MalyshevMinlos2}).

Unfortunately, eigenvectors ${{\psi}''_{\alpha}}$
of the operators $H''_{\alpha}$  differ from those
for the initial spectral problem~(\ref{ini}):
${{\psi}''}_{\alpha}=X_{\alpha}^{1/2}\psi_{\alpha}$.

\begin{lemma}\label{LHolder}
 Let the kernel $B_{\beta\alpha}(\mu,\lambda)$, 
${\beta}\neq{\alpha}$,
of the operator $B_{\beta\alpha}$ belong 
to the class ${\cal B}^{{\beta}{\alpha}}_{\theta\gamma}$
with $\theta > \frac{1}{2}$
and $Q_{\beta\alpha}$ be a solution of 
Eq.~(\ref{Qbasic1}) satisfying together
with $Q_{\alpha\beta}=-Q_{\beta\alpha}^{*}$ the conditions
(\ref{resolvest}). Then

(a) the operator $Q_{\beta\alpha}$ is an integral operator, 
$Q_{\beta\alpha}:$
${\cal H}_{\alpha}\rightarrow{\cal H}_{\beta},$ 
with a kernel $Q_{\beta\alpha}(\mu,\lambda)$
belonging to ${\cal B}^{{\beta}{\alpha}}_{\theta\gamma}$;

(b) the potential 
$W_{\alpha}\equiv B_{\alpha\beta} Q_{\beta\alpha}$ is
an integral operator, $W_{\alpha}:$
${\cal H}_{\alpha}\rightarrow{\cal H}_{\alpha},$ 
with a kernel $W_{\alpha}(\lambda,\lambda')$
belonging to ${\cal B}^{{\alpha}{\alpha}}_{\theta\gamma}$.
\end{lemma}

\noindent{\sl Proof.} At the beginning we prove the assertion (b).
According to~(\ref{Qadj}),
\be
\label{w}
W_{\alpha}=-B_{\alpha\beta}\int\limits_{\sigma_{\alpha}}
(H_{\beta}^{*}-\lambda)^{-1}B_{\beta\alpha} E_{\alpha}(d\lambda)
\ee
with $H_{\beta\alpha}^{*}=A_{\beta}+W_{\beta}^{*}=
A_{\beta}+Q_{\alpha\beta}^{*}B_{\alpha\beta}$. Since the inequalities
(\ref{resolvest}) take place we write
$$
\| (H_{\beta}^{*}-\lambda)^{-1} \| = \| (H_{\beta}-\lambda)^{-1} \|
\leq C_{\beta\alpha}
$$
for any $\lambda\in\sigma_{\alpha}$. In the diagonal representation
(\ref{Neumann}),(\ref{multi}), the equation~(\ref{w}) turns in
$$
W_{\alpha}(\lambda,\lambda')=
-B_{\alpha\beta}(\lambda,\,\cdot\,)(H_{\beta}^{*}-\lambda)^{-1}
B_{\beta\alpha}(\,\cdot\, ,\lambda').
$$
It means that
$$
|W_{\alpha}(\lambda,\lambda')| \leq
\|B_{\alpha\beta}(\lambda,\,\cdot\, )\|_{{\cal H}_{\beta}}\cdot
\|(H_{\beta}^{*}-\lambda)^{-1}\|\cdot
\| B_{\beta\alpha}(\,\cdot\, ,\lambda')\|_{{\cal H}_{\beta}}\leq
$$
\be
\label{West}
\leq C_{\beta\alpha} 
\|B_{\alpha\beta}(\lambda,\,\cdot\, )\|_{{\cal H}_{\beta}}\cdot
\| B_{\beta\alpha}(\,\cdot\, ,\lambda')\|_{{\cal H}_{\beta}}.
\ee
Here,
$
\|B_{\alpha\beta}(\lambda,\,\cdot\, )\|_{{\cal H}_{\beta}}=
\left[  \mathop{\mbox{$\sum$\hskip-1em{\LARGE $\int$}\hskip0.1em}}
\limits_{\sigma_{\beta}}
|B_{\alpha\beta}(\lambda,\mu)|^{2}d\mu \right]^{1/2}.
$
Since $\theta >\frac{1}{2},$ we have
$
\|B_{\alpha\beta}(\lambda,\,\cdot\, )\|_{{\cal H}_{\beta}}
$ $\leq$ $ \frac{c(\theta)}{(1+|\lambda|)^\theta}\cdot \|B\|_{{\cal B}}
$
with some $c(\theta)$, $\, c(\theta) > 0$, depending only on $\theta$.
Analogously,
$$
\|B_{\beta\alpha}(\,\cdot\, ,\lambda')\|_{{\cal H}_{\beta}}=
\|\overline{B_{\alpha\beta}}(\lambda',\,\cdot\,)\|_{{\cal H}_{\beta}}
\leq \frac{c(\theta)}{(1+|\lambda|)^\theta}\cdot \|B\|_{{\cal B}},
$$
where the operator
$\overline{B_{\alpha\beta}}(\lambda,\mu),$ 
$\overline{B_{\alpha\beta}}(\lambda,\mu):$
${\cal G}_{\beta}(\mu)\rightarrow{\cal G}_{\alpha}(\lambda)$,
is adjoint to $B_{\beta\alpha}(\mu,\lambda)$.

Estimations similar to~(\ref{West}) may be done also for
$|W_{\alpha}(\lambda'',\lambda')-W_{\alpha}(\lambda,\lambda')|,$
$\, \lambda,\lambda''\in\sigma_{\alpha}^c,$
$\, \lambda\in\sigma_{\alpha}$,
$|W_{\alpha}(\lambda,\lambda''')-W_{\alpha}(\lambda,\lambda')|,$
$\, \lambda\in\sigma_{\alpha}$,
$\, \lambda''',\lambda'\in\sigma_{\alpha}^c,$
and
$ |W_{\alpha}(\lambda,\lambda')-W_{\alpha}(\lambda'',\lambda') -$
$W_{\alpha}(\lambda,\lambda''')+W_{\alpha}(\lambda'',\lambda''')|,\,$
$\, \lambda,\lambda',\lambda'',\lambda'''\in\sigma_{\alpha}^c$,
in terms of the norms
$
\|B_{\alpha\beta}(\lambda,\,\cdot\,) - 
B_{\alpha\beta}(\lambda'',\,\cdot\,)\|_{{\cal H}_{\beta}}
$
                and
$
\|B_{\beta\alpha}(\,\cdot\, ,\lambda''') - 
B_{\beta\alpha}(\,\cdot\, ,\lambda')\|_{{\cal H}_{\beta}}.
$
Estimating the latter through $ \| B_{\alpha\beta}\|_{\cal B} $
we come to the inequality
$$
  \| W_{\alpha}\|_{{\cal B}^{{\alpha}{\alpha}}_{\theta\gamma}}
  \leq c(\theta)C_{\beta\alpha}\cdot
  \| B_{\beta\alpha}\|^2_{{\cal B}^{{\beta}{\alpha}}_{\theta\gamma}}
$$
with  $0< c(\theta) < \infty. $  
Therefore, we have proved the assertion (b).

To prove the statement (a) we note that according to~(\ref{Qbasic1}),
$$
Q_{\beta\alpha}=
\int\limits_{\sigma_{\beta}} E_{\beta}(d\mu) B_{\beta\alpha} 
\left[
(A_{\alpha} -\mu)^{-1}-(H_{\alpha}-\mu)^{-1}
W_{\alpha}(A_{\alpha}-\mu)^{-1}  \right]
$$
or, in the diagonal representation~(\ref{Neumann}),(\ref{multi}),
$$
Q_{\beta\alpha}(\mu,\lambda)=
\frac{B_{\beta\alpha}(\mu,\lambda)}{\lambda-\mu}
- \frac{B_{\beta\alpha}(\mu,\,\cdot\,)
(H_{\alpha}-\mu)^{-1}W_{\alpha}(\,\cdot\, ,\lambda)}
{\lambda-\mu}.
$$
Repeating literally the last part of the proof of the assertion (b)
we come to the inequality
$$
\| Q_{\beta\alpha}\|_{{\cal B}^{{\beta}{\alpha}}_{\theta\gamma}} \leq
\mathop{\rm sup}\limits_{
                           \mbox{
                                  \scriptsize
                                           $\begin{array}{c}
                                                 \mu\in\sigma_{\beta}\\
                                             \lambda\in\sigma_{\alpha}
                                             \end{array}$
                                 }
    }
     \frac{1}{   |\lambda-\mu| } \cdot
                                                 \left\{
 \| B_{\beta\alpha} 
\|_{  {\cal B}^{{\beta}{\alpha}}_{\theta\gamma}   } +
 \mbox{\phantom{
$\mathop{\rm sup}\limits_{\mu\in\sigma_{\beta}}$ }}\right.
$$
$$
    \left. \mbox{\phantom{{\Large I}}}
+ c(\theta)\cdot
\| B_{\beta\alpha}\|_{{\cal B}^{{\beta}{\alpha}}_{\theta\gamma}} 
\cdot
\mathop{\rm sup}\limits_{\mu\in\sigma_{\beta}}  
\| (H_{\alpha}-\mu)^{-1}  \|\cdot
\| W_{\alpha}\|_{{\cal B}^{{\alpha}{\alpha}}_{\theta\gamma}}
 \right\},
               \quad 0 < c(\theta) < +\infty.
$$
Consequently 
$Q_{\beta\alpha}\in{\cal B}^{{\beta}{\alpha}}_{\theta\gamma} $
and
$$
\|Q_{\beta\alpha}\|_{{\cal B}} \leq \frac{1}{d_0}\cdot
       \left\{
\|B_{\beta\alpha}\|_{{\cal B}}  +
c(\theta)C_{\alpha\beta} C_{\beta\alpha}
\cdot \|B_{\beta\alpha}\|_{{\cal B}}^3                
               \right\},
\quad 0 <c(\theta) < +\infty.
$$
This completes the proof of Lemma~\ref{LHolder}.

\noindent\begin{corollary}\label{CHolderSmoothness}
{\it If $B_{\beta\alpha}\in{\cal B}^{{\beta}{\alpha}}_{\theta\gamma}$,
$\,\theta > \frac{1}{2},$
then the solution of Eq.~(\ref{Qbasic1}) described by
Theorem~\ref{ThSolvability}
belongs to the class 
${\cal B}^{{\beta}{\alpha}}_{\theta\gamma}$, too.   }
\end{corollary}

This statement is based on the fact that the mentioned solution satisfies
automatically the conditions~(\ref{BQd}) and, hence, the conditions
(\ref{resolvest}).
\section{EIGENFUNCTIONS AND THE EXPANSION THEOREM}
\label{Eigenfunctions}
In the preceding section, we have proved the existence
(in the unit ball 
$M_{\beta\alpha}(1)\subset
{\bf B}({\cal H}_{\alpha},{\cal H}_{\beta})$)
of a solution
$Q_{\beta\alpha}$ of the basic 
equation~(\ref{Qbasic1}) only in the case where spectra
$\sigma_1,$ $\sigma_2$ of the operators $A_1$, $A_2$ are separated,
${\rm dist}\{  \sigma_1,\sigma_2\}$ = $d_0$ $> 0$, and
$\|B_{12}\|_2 = \| B_{21} \|_2 < \frac{d_0}{2}$. May be, however,
Eqs.~(\ref{basic}) and~(\ref{Qbasic1}) have solutions also in other
cases. That is why we study the spectral properties of the operator
$H_{\alpha}=A_{\alpha}+B_{\alpha\beta} Q_{\beta\alpha}$ 
not supposing that
$\|B_{\alpha\beta}\|_2 < \frac{d_0}{2}$
and using more general requirements~(\ref{resolvest}) only, with
$C_{\alpha\beta}$,  some positive numbers, 
${\alpha},{\beta}=1,2$, ${\beta}\neq{\alpha}$. Of course, we
assume again that the condition~(\ref{dist}) takes place.
Remember that the requirements~(\ref{resolvest}) are sufficient for
existence of the operators 
$V_{\alpha}(H_{\alpha})$. As well, the equations
(\ref{QbasicSym}) and~(\ref{Qadj}) take place and the assertions of
Theorem~\ref{ThInvariant} and Lemma~\ref{LHolder} are valid.

So, let us suppose that 
$Q_{\beta\alpha}$ and $Q_{\alpha\beta}=-Q_{\beta\alpha}^{*}$ 
are solutions
of Eqs.~(\ref{Qbasic1}) and~(\ref{QbasicSym}) satisfying the
conditions~(\ref{resolvest}). 
It follows from Lemma~\ref{LQdefinition} that
$Q_{\beta\alpha}\in{\bf B}_{\beta\alpha}
({\cal H}_{\alpha},{\cal H}_{\beta})$ 
as well as 
$Q_{\alpha\beta}\in{\bf B}_{\alpha\beta}
({\cal H}_{\beta},{\cal H}_{\alpha})$.
If $B_{\beta\alpha}\in{\cal B}_{\theta\gamma}^{{\beta}{\alpha}}$, 
$\theta>\frac{1}{2}$, then,
according to Lemma~\ref{LHolder}, 
$Q_{\beta\alpha}\in{\cal B}_{\theta\gamma}^{{\beta}{\alpha}}$
and $Q_{\alpha\beta}\in{\cal B}_{\theta\gamma}^{{\alpha}{\beta}}$.

By Theorem~\ref{ThInvariant}, the operator
${{\bf H}}'=\mathop{\rm diag}\{ H_1, H_2  \}$ is connected with
the (self-adjoint) operator  ${{\bf H}}$
by a similarity transform. Thus, the spectra
$\sigma(H_1)$ and $\sigma(H_2)$ of 
the operators $H_{\alpha}$, ${\alpha}=1,2$,
are real and 
$\sigma(H_1)\bigcup\sigma(H_2)=\sigma({{\bf H}})$. 
Continuous spectrum
$\sigma_c (H_{\alpha})$ 
of the each operator $H_{\alpha}$ coincides with that of the
operator $A_{\alpha}$, $\sigma_c (H_{\alpha})=\sigma_{\alpha}^c$,
since due to $\|B_{\alpha\beta}\|_2 < +\infty$, 
the potential $W_{\alpha}=B_{\alpha\beta} Q_{\beta\alpha} $ is
a compact operator. Since $\sigma_1^c \bigcap\sigma_2^c =\emptyset $
we have $\sigma_c (H_1)\bigcap\sigma_c (H_2) =\emptyset$. 
We show now that the
discrete spectra $\sigma_d (H_{\alpha}),$ ${\alpha}=1,2$, 
satisfy a similar condition.

Let us suppose that $\sigma _{d}(H_{\alpha}  )\neq\emptyset $,
$z\in \sigma _{d}(H_{\alpha}  )$ and $\psi _{\alpha}  $  is
the  corresponding  eigenfunction  of $H_{\alpha}  ,$
$ H_{\alpha}  \psi _{\alpha}  =z\psi _{\alpha}  ,$
$\psi _{\alpha}  \in {\cal D}(H_{\alpha}  )={\cal D}(A_{\alpha}  )$.
Then, according to construction  of $H_{\alpha}$,
we   have
$H_{\alpha}  \psi _{\alpha}  =$
$(A_{\alpha}  +V_{\alpha}  (H_{\alpha}  ))\psi _{\alpha}  =$
$(A_{\alpha}  +V_{\alpha}  (z))\psi _{\alpha}  =z\psi _{\alpha}  $.
Thus if $z\in \sigma _{d}(H_{\alpha}  )$  then $z$  becomes
automatically a point  of the discrete
spectrum of the  initial spectral problem ~(\ref{ini}).
At the same time $\psi_{\alpha} $ becomes it's eigenfunction.

Let us further denote the eigenfunctions  of  the
operator $H_{\alpha}  $  
discrete  spectrum  by $\psi ^{(j)} _{\alpha}  , $
$\psi ^{(j)} _{\alpha}  =u^{(j)} _{\alpha}  $,
keeping for them the  same  numeration  as  for
eigenvectors of $U_j$, $U_j=\{ u_{\alpha}^{(j)}, u_{\beta}^{(j)} \}, $
of the Hamiltonian ${{\bf H}}$,
${{\bf H}}U_j =z_j U_j$, $z_j\in\sigma_d ({{\bf H}}).$ 
We assume that in the case
of multiple discrete eigenvalues, certain $z_j$ may be repeated
in this numeration.
 By ${\cal U} ^{d}$
we denote the set ${\cal U} ^{d}=\{U_{j}, j=1,2,\ldots\}$
of  all the eigenvectors $U_j$.

Let ${\cal U} ^{d}_{\alpha}  $
be such a subset of ${\cal U} ^{d}$
that it's elements have
the operator $H_{\alpha}$ eigenvectors $\psi^{(j)}_{\alpha}$
in the  capacity  of the  channel ${\alpha}$
components:
${\cal U} _{\alpha}^d =$ $\{U_{j}:$
$ U_{j}=\{u^{(j)} _{1},u^{(j)} _{2}\},$
$ u^{(j)} _{\alpha}  =\psi ^{(j)} _{\alpha}  \}$.
By Theorem~\ref{ThInvariant}, we have
${\cal U} ^{d}_1  \bigcup {\cal U} ^{d}_2 = {\cal U} ^{d}.$

\begin{theorem}\label{ThDiscSpectrum}
Let $H_{\beta}  =A_{\beta}  + B_{\beta\alpha} Q_{\alpha\beta} ,$
correspond  (for
$\| B_{\beta\alpha} \|_{2} < +\infty $)
to the same  solution $Q_{\alpha\beta} =- Q_{\beta\alpha}^{*} $
of  Eqs.~(\ref{Qbasic1}) and~(\ref{QbasicSym})  as
$H_{\alpha} =A_{\alpha} +B_{\alpha\beta} Q_{\beta\alpha} $,
and the conditions~(\ref{resolvest}) are valid.
Let $z_{j}\in \sigma _{d}(H_{\alpha} )$ and 
$H_{\alpha} u^{(j)} _{\alpha} =
z_{j}u^{(j)} _{\alpha} $ with $u_{\alpha}^{(j)} $,
the channel ${\alpha}$ component of the eigenvector
$U_j =\{ u_{\alpha}^{(j)},u_{\beta}^{(j)} \}$ 
of the operator ${{\bf H}}$, ${{\bf H}}U_j =z_j U_j$.
Then either
 $z_{j}\not\in \sigma _{d}(H_{\beta})$, ${\beta}\neq{\alpha}$,
 or  (if $z_{j}\in\sigma_{d}(H_{\beta}  ))$ 
the vector $u_{\beta}^{(j)} $ is not
an eigenvector of $H_{\beta}$.
\end{theorem}

\begin{corollary}\label{CDiscSpectrum}
$
{\cal U} ^{d}_1
\bigcap{\cal U} ^{d}_2 =\emptyset
$.
\end{corollary}

Statement of Theorem~\ref{ThDiscSpectrum} means that discrete
spectrum $\sigma_d ({{\bf H}}) $ is distributed between
discrete spectra $\sigma_d (H_1)$ and 
$\sigma_d (H_2)$ in such a way that
operators $H_1$ and $H_2$ have not ``common'' eigenvectors
$U_j =\{ u_1^{(j)},u_2^{(j)} \}$: 
simultaneously, component $u_1^{(j)}$
can not be eigenvector
for $H_1$, and  $u_2^{(j)}$ with the same $j$,  for $H_2$.
\medskip

\noindent{\sl Proof} of the Theorem will be given by contradiction.

Let us suppose that $\psi_{\alpha}^{(j)} = u_{\alpha}^{(j)}$ 
is an eigenvector of $H_{\alpha}$
corresponding to $z_j$ i.e.
\be
\label{Aeigen}
(A_{\alpha}+B_{\alpha\beta} Q_{\beta\alpha} -z_j)
\psi_{\alpha}^{(j)}=0.
\ee

If $z_j\in\sigma_{\alpha}=\sigma(A_{\alpha})$ then automatically
$z_j\not\in\sigma_d(H_{\beta})$ 
since due to conditions~(\ref{resolvest})
we have $\sigma(H_{\beta})\bigcap\sigma(A_{\alpha})=\emptyset$. 
Thus in the case
where $z_j\in\sigma_{\alpha}$ the assertion of Theorem is valid.

Let $z_j\not\in\sigma(A_{\alpha})$. In this case we can rewrite
Eq.~(\ref{Aeigen}) in the form
\be
\label{LSch}
\psi_{\alpha}^{(j)} = 
-(A_{\alpha} -z_j)^{-1}B_{\alpha\beta} 
Q_{\beta\alpha}\psi_{\alpha}^{(j)}.
\ee
Let $y_{\beta}^{(j)}=Q_{\beta\alpha}\psi_{\alpha}^{(j)}$. 
It follows from~(\ref{LSch}) that
\be
\label{y}
y_{\beta}^{(j)} + 
Q_{\beta\alpha}(A_{\alpha} - z_j)^{-1} y_{\beta}^{(j)} =0.
\ee
We will show that the vector $ y_{\beta}^{(j)} $ 
is a solution of the initial
spectral problem~(\ref{ini}) in the channel ${\beta}$ at $z=z_j$ and
$\tilde{U}_j=\{ \psi_{\alpha}^{(j)}, y_{\beta}^{(j)} \}$ 
is an eigenvector of ${{\bf H}}$,
${\bf H}\tilde{U}_j=z_j \tilde{U}_j$. 
To do this, we act on both parts of
Eq.~(\ref{y}) by    $ H^{*}_{\beta}-z_j$  
remembering that, according to
(\ref{Qadj}),
$Q_{\beta\alpha}=-Q_{\alpha\beta}^{*}$
$=-\int\limits_{\sigma_{\alpha}}(H_{\beta}^{*}-\lambda)^{-1}
B_{\beta\alpha} E_{\alpha}(d\lambda)$.
We obtain
$$
(H_{\beta}^{*}-z_j)y_{\beta}^{(j)}+
\int\limits_{\sigma_{\alpha}} (H_{\beta}^{*}-z_j)
(H_{\beta}^{*}-\lambda)^{-1}(z_j -\lambda)^{-1} 
B_{\beta\alpha} E_{\alpha}(d\lambda)
   B_{\alpha\beta} y_{\beta}^{(j)} =0.
$$
Using the identity
$ (H-z)(H-\lambda)^{-1}(z-\lambda)^{-1}=
  (z-\lambda)^{-1}-(H-\lambda)^{-1} $
we find
$$
 (H_{\beta}^{*}-z_j)y_{\beta}^{(j)}+\int\limits_{\sigma_{\alpha}}
 [(z_j -\lambda)^{-1} - (H_{\beta}^{*}-z_j) ]
 B_{\beta\alpha} E_{\alpha}(d\lambda) 
B_{\alpha\beta} y_{\beta}^{(j)} =0
$$
or, and it is the same,
\be
\label{yy}
  (H_{\beta}^{*}-z_j)y_{\beta}^{(j)} - 
B_{\alpha\beta}(A_{\alpha} -z_j)^{-1} B_{\alpha\beta} y_{\beta}^{(j)}
   + Q_{\beta\alpha} B_{\alpha\beta} y_{\beta}^{(j)} =0.
\ee
However 
$ H_{\beta}^{*}=A_{\beta} -
Q_{\beta\alpha} B_{\alpha\beta} $. 
Hence the relation~(\ref{yy})
turns in equation~(\ref{ini}) for the channel ${\beta}$,
$$
  [A_{\beta} - B_{\beta\alpha} 
(A_{\alpha} -z_j)^{-1} B_{\alpha\beta} - z_j ] y_{\beta}^{(j)} = 0.
$$
So, we have proved that $y_{\beta}^{(j)}$ 
is a solution of the initial problem
in the channel ${\beta}$ and we did deal with an eigenvector
$ U_j=\{ u_{\alpha}^{(j)}, u_{\beta}^{(j)} \} $ 
of the operator ${{\bf H}}$ having the components
$u_{\alpha}^{(j)}=\psi_{\alpha}^{(j)}$ and 
$u_{\beta}^{(j)} =y_{\alpha}^{(j)} $.

Let us show that $y_{\beta}^{(j)}$ can not be 
an eigenvector of $H_{\beta}$
corresponding to the eigenvalue $z_j$. 
Actually, due to~(\ref{y}) we have
$$
{\rm a}\equiv \langle y_{\beta}^{(j)} +
Q_{\beta\alpha} (A_{\alpha} -z_j)^{-1}
B_{\alpha\beta} y_{\beta}^{(j)},
 y_{\beta}^{(j)}  \rangle  =0.
$$
On the other hand
$$
{\rm a}= \| y_{\beta}^{(j)}\|^2 + 
\langle (A_{\alpha}-z_j)^{-1} B_{\alpha\beta} y_{\beta}^{(j)},
  Q_{\beta\alpha}^{*}y_{\beta}^{(j)}   \rangle.
$$
If $y_{\beta}^{(j)}$ is an eigenvector of 
$H_{\beta}$, $H_{\beta} y_{\beta}^{(j)} =z_j y_{\beta}^{(j)}$,
then
$$
Q_{\beta\alpha}^{*} y_{\beta}^{(j)} =
-Q_{\alpha\beta} y_{\beta}^{(j)} = 
-\int\limits_{\sigma_{\alpha}}
E_{\alpha}(d\lambda) B_{\alpha\beta} 
(H_{\beta} -\lambda)^{-1} y_{\beta}^{(j)} =
  (A_{\alpha} -z_j)^{-1} B_{\alpha\beta} y_{\beta}^{(j)}.
$$
It means that
$$
{\rm a} = \|y_{\beta}^{(j)}\|^2 +
\| (A_{\alpha}-z_j)^{-1} B_{\alpha\beta} y_{\beta}^{(j)}\|^2
\geq  \|y_{\beta}^{(j)}\|^2 .
$$
Since ${\rm a} =0 $ we get  $y_{\beta}^{(j)} =0 $
and, due to~(\ref{LSch}), $\psi_{\alpha}^{(j)} =0.$ 
However, by supposition,
$ \psi_{\alpha}^{(j)} \neq 0.$
 Thus, we come to a contradiction and 
$y_{\beta}^{(j)}$ can not be
an eigenvector of $H_{\beta}$. 
And so, if $z_j \in \sigma_d (H_{\alpha})$ and
$H_{\alpha} u_{\alpha}^{(j)} =z_j u_{\alpha}^{(j)}$ 
then $u_{\beta}^{(j)} $ is not an eigenvector of $H_{\beta}$.
The proof of Theorem~\ref{ThDiscSpectrum} is completed.
\medskip

Let us pay attention to the continuous spectrum of $H_{\alpha}$
assuming here that 
$B_{\alpha\beta}\in{\cal B}^{{\beta}{\alpha}}_{\theta\gamma},$
$\theta > \frac{1}{2},$ $\gamma > \frac{1}{2}$, and consequently,
$Q_{\alpha\beta}\in{\cal B}^{{\beta}{\alpha}}_{\theta\gamma},$ 
${\alpha},{\beta}=1,2$, ${\beta}\neq{\alpha}.$

Consider at $\lambda'\in\sigma_{\alpha}^c$ 
the integral equations
\be
\label{LSchpsi}
\psi ^{(\pm )}_{\alpha}  (\lambda ,\lambda ')=
I_{\alpha}^c  {\delta}  (\lambda -\lambda ')-
[(A_{\alpha}  -\lambda '\mp i0)^{-1}
W_{\alpha}
\psi ^{(\pm )}_{\alpha}](\lambda ,\lambda '), \,\, {\alpha}=1,2,
\ee
where as usually 
$W_{\alpha}=B_{\alpha\beta} Q_{\beta\alpha}.$ 
Since
$W_{\alpha}\in{\cal B}^{{\alpha}{\alpha}}_{\theta\gamma},$
the integral operator with the kernel
$
\frac{W_{\alpha}(\lambda,\lambda')}{\lambda -\lambda '\mp i0}
$
is compact  in
$
{\cal M}_{\theta '\gamma '},
$
$
\frac{1}{2}<\theta '<\theta ,
$
$
0<\gamma '<\gamma
$
(cf. Refs.~\onlinecite{LadyzhFaddeev}, \onlinecite{Faddeev64}).
If $\lambda '\not\in \sigma _{d}(H_{\alpha}  )$
then Eq.~(\ref{LSchpsi}) for
$\psi^{(+)}_{{\alpha}}$
as well as for $\psi^{(-)}_{{\alpha}}$ is
uniquely  solvable (see Ref.~\onlinecite{Faddeev64}) in
the  class  of  the  form ~(\ref{amplitude})
distributions.

Denote by
$
\Psi ^{(\pm )}_{\alpha}  ,
$
$\Psi ^{(\pm )}_{\alpha}  :$
$
{{\cal H}} ^{c}_{\alpha}  \rightarrow {{\cal H}} _{\alpha}
$,
the in\-te\-gral ope\-ra\-tor
with the  ker\-nel
$
\psi ^{(\pm )}_{\alpha}  (\lambda ,\lambda ')
$.
The  operator $\Psi ^{(\pm )}_{\alpha}  $  is
bounded and
$
\Psi ^{(\pm )}_{\alpha}  {\cal D}(A^{(0)}_{\alpha}  )
\subseteq {\cal D}(H_{\alpha}  )$%
~\cite{LadyzhFaddeev}, \cite{Faddeev64}.
It  follows  from ~(\ref{LSchpsi})
that $\Psi ^{(\pm )}_{\alpha}  $
has   the   property
$
H_{\alpha}  \Psi ^{(\pm )}_{\alpha}  =
\Psi ^{(\pm )}_{\alpha}  A^{(0)}_{\alpha}
$.
Thus,
$
Q_{\beta\alpha} \Psi ^{(\pm )}_{\alpha}  (\cdot ,\lambda ')=
$
$
(\lambda '- A_{\beta} )^{-1}
$
$
B_{\beta\alpha} \Psi ^{(\pm )}_{\alpha}  (\cdot ,\lambda ')
$.
Substitution    of    this
expression in~(\ref{LSchpsi}) shows that
$
\psi ^{(\pm )}_{\alpha}
$
satisfies~(\ref{wavefunction}).
Due to  the  uniqueness
of  Eq.~(\ref{wavefunction}) solution at 
$\lambda'\not\in\sigma_d({{\bf H}})$
  we   have
$
\psi ^{(\pm )}_{\alpha}  (\lambda ,\lambda ')=
u^{(\pm )}_{{\alpha}{\alpha}}(\lambda,\lambda')
$.
This    means    that    each
eigenfunction
$
u^{(\pm )}_{{\alpha} {\alpha} }(\lambda ,\lambda '),
$
$
\lambda ^\prime \in \sigma ^{c}_{\alpha}  ,
$
$
\lambda '\not\in \sigma _{d}({{\bf H}})
$
of the initial spectral problem
problem~(\ref{ini}) is also an eigenfunction  of $H_{\alpha}  $.

Consider  the  functions
$
\tilde{\psi}^{(j)} _{\alpha}  =
\psi ^{(j)} _{\alpha}  -Q_{\alpha\beta} u^{(j)} _{\beta}
$
and \newline
$
\tilde{\psi}^{(\pm )}_{\alpha}  (\cdot ,\lambda ')=
$
$
\psi ^{(\pm )}_{\alpha}  (\,\cdot\, ,\lambda ')-
$
$
Q_{\alpha\beta} u^{(\pm )}_{\beta\alpha}
(\,\cdot\, ,\lambda '),
$
$
\lambda '\in \sigma ^{c}_{\alpha}
$.
Let
$
\tilde{\Psi }^{(\pm )}_{\alpha}  ,
$
$
\tilde{\Psi }^{(\pm )}_{\alpha}  :
$
$
{{\cal H}} ^{c}_{\alpha}  \rightarrow {{\cal H}} _{\alpha}
$,
be the integral operator  with  the
kernel
$
\tilde{\psi}^{(\pm )}_{\alpha}  (\lambda ,\lambda ')
$.

\begin{theorem}\label{ThExpansion}
 The functions
$
\tilde{\psi}^{(j)} _{\alpha}
$
(with $j$  such  that
$
U_{j}\in {\cal U} ^{d}_{\alpha}  )
$
 are eigenfunctions of adjoint  operator
$
H^{*}_{\alpha}  , H^{*}_{\alpha}  =
A_{\alpha}  +Q^{*}_{\beta\alpha} B_{\beta\alpha} ,
$
  discrete  spectrum,
$
H^{*}_{\alpha}  \tilde{\psi}^{(j)} _{\alpha}  =
$
$
z_{j}\tilde{\psi}^{(j)} _{\alpha}
$.
 Operators
$
\tilde{\Psi }^{(\pm )}_{\alpha}
$
 have the property
$
H^{*}_{\alpha}  \tilde{\Psi }^{(\pm )}_{\alpha}=
$
$
\tilde{\Psi }^{(\pm )}_{\alpha}  A^{(0)}_{\alpha}
$.
 At the same  time  the  orthogonality  relations  take
place:
$
\langle\psi^{(j)} _{\alpha}  ,\tilde{\psi}^{(k)}_{\alpha}  \rangle=
{\delta}  _{jk},
$
$
\Psi ^{(\pm )*}_{\alpha}  \tilde{\Psi }^{(\pm )}_{\alpha}  =
\left.{I_{\alpha}}\right|_{{\cal H}_{\alpha}^c},
$
$
\tilde{\Psi }^{(\pm )*}_{\alpha}  \psi ^{(j)} _{\alpha}  =0
$
                       and
$
\Psi_{\alpha}^{(\pm)*}\tilde{\psi}^{(j)} _{\alpha}  =0.
$
  Also,  the following  completeness  relations are valid,
\be
\label{completeness}
\sum\limits_{j:\, U_{j}\in{\cal U} _{\alpha}^d}
{\psi}^{(j)} _{\alpha} \langle\cdot,
\tilde{\psi}^{(j)} _{\alpha}  \rangle+
\Psi^{(\pm )}_{\alpha}  \tilde{\Psi }^{(\pm )*}_{\alpha}  =
I_{\alpha}  , \quad {\alpha} =1,2,
\ee
\end{theorem}

\noindent{\sl Proof.} Show for example that
\be
\label{Hadjeigen}
  H_{\alpha}^{*}\tilde{\psi}_{\alpha}^{(j)}=
z_j\tilde{\psi}_{\alpha}^{(j)}
\ee
(remember that $z_j\in{\bf R}$). We have
$$
H_{\alpha}^{*}\tilde{\psi}_{\alpha}^{(j)}=
(A_{\alpha}-Q_{\alpha\beta} B_{\beta\alpha})(\tilde{\psi}_{\alpha}^{(j)} 
-Q_{\alpha\beta} u_{\beta}^{(j)})=
$$
\be
\label{Heq1}
=(A_{\alpha}-Q_{\alpha\beta} 
B_{\beta\alpha})\psi_{\alpha}^{(j)} -
(A_{\alpha} Q_{\alpha\beta}-
Q_{\alpha\beta} B_{\beta\alpha} Q_{\alpha\beta})u_{\beta}^{(j)}.
\ee
Note that $A_{\alpha}=H_{\alpha}-
B_{\alpha\beta} Q_{\beta\alpha}$ and, hence,
\be
\label{Heq2}
(A_{\alpha}-Q_{\alpha\beta} 
B_{\beta\alpha})\psi_{\alpha}^{(j)} =
z_j\psi_{\alpha}^{(j)} -
(B_{\alpha\beta} Q_{\beta\alpha}+
Q_{\alpha\beta} B_{\beta\alpha})\psi_{\alpha}^{(j)}.
\ee
The second term in the right part of~(\ref{Heq2}) may be easily
expressed by $u_{\beta}^{(j)}$. Actually,
$
u_{\beta}^{(j)}=-(A_{\beta}-z_j)^{-1}B_{\beta\alpha} 
\psi_{\alpha}^{(j)} $
[here, we use again the
property $\sigma(H_{\alpha})\bigcap\sigma_{\beta}=\emptyset$ 
following
from~(\ref{resolvest})]. Since Eqs.~(\ref{Qbasic}) and
$H_{\alpha}\psi_{\alpha}^{(j)}=z_j\psi_{\alpha}^{(j)}$ 
take place, we find
$Q_{\beta\alpha}\psi_{\alpha}^{(j)}=
B_{\alpha\beta} u_{\beta}^{(j)}.$ Thus,
$$
 (B_{\alpha\beta} Q_{\beta\alpha} +
Q_{\alpha\beta} B_{\beta\alpha})\psi_{\alpha}^{(j)}=
    B_{\beta\alpha} Q_{\beta\alpha}\psi_{\alpha}^{(j)}+
Q_{\alpha\beta}(A_{\beta}-z_j)
(A_{\beta}-z_j)^{-1}B_{\beta\alpha}\psi_{\alpha}^{(j)}=
$$
$$
   =B_{\beta\alpha} u_{\beta}^{(j)} -Q_{\alpha\beta} 
(A_{\beta} - z_j)u_{\beta}^{(j)}.
$$
Substituting the expressions obtained into~(\ref{Heq2}) and then
into~(\ref{Heq1}), we get
$$
  H_{\alpha}^{*}\tilde{\psi}_{\alpha}^{(j)}=
z_j(\psi_{\alpha}^{(j)}-Q_{\alpha\beta} u_{\beta}^{(j)})
 + [-B_{\alpha\beta}+Q_{\alpha\beta} A_{\beta}- 
A_{\alpha} Q_{\alpha\beta} +
Q_{\alpha\beta} B_{\beta\alpha} Q_{\alpha\beta}]u_{\beta}^{(j)}.
$$
According to the equations~(\ref{QbasicSym}), the expression in the square
brackets is equal to zero and we come to~(\ref{Hadjeigen}).

The equalities
$
     H_{\alpha}^{*}\tilde{\psi}^{(\pm)}_{\alpha}(\,\cdot\, ,\lambda')=
$
$
     \lambda'\tilde{\psi}^{(\pm)}_{\alpha}(\,\cdot\, ,\lambda'),
$
$\lambda'\in\sigma_{\alpha}^c ,$ are proved quite analogously.

The orthogonality relations
$
\langle\psi^{(j)} _{\alpha}  ,\tilde{\psi}^{(k)}_{\alpha}  \rangle=
{\delta}_{jk},
$
$
\tilde{\Psi }^{(\pm )*}_{\alpha}  \psi ^{(j)} _{\alpha}  =0
$
                       and
$
\Psi_{\alpha}^{(\pm)*}\tilde{\psi}^{(j)} _{\alpha}  =0
$
are trivial. Proofs of the relation
$
\Psi ^{(\pm )*}_{\alpha}  \tilde{\Psi }^{(\pm )}_{\alpha}  =
\left.{I_{\alpha}}\right|_{{\cal H}_{\alpha}^c},
$
and the equality~(\ref{completeness}) are very similar. 
Both these proofs are
based on use of properties of the wave 
operators $U^{(\pm)}$. As a sample,
we give a proof of the completeness relation~(\ref{completeness}).

Consider the operator
$$
{\cal A}=\sum\limits_{j:\, U_j\in{\cal U}_{\alpha}^d}
\psi_{\alpha}^{(j)} \langle\,\cdot\, ,
\tilde{\psi}_{\alpha}^{(j)}\rangle +
\Psi_{\alpha}\tilde{\Psi}_{\alpha}^{*} =
$$
\be
\label{Aeq0}
=\sum\limits_{j:\, U_j\in{\cal U}_{\alpha}^d}
\psi_{\alpha}^{(j)} \langle\,\cdot\, ,
\psi_{\alpha}^{(j)} -Q_{\alpha\beta} u_{\beta}^{(j)}\rangle +
\Psi_{\alpha}[\Psi_{\alpha}^{*} -(Q_{\alpha\beta} u_{\beta\alpha})^{*}].
\ee
For convenience, we omit signs ``$\pm$'' in notations of
${\Psi}^{(\pm)}_{\alpha} \equiv  u^{(\pm)}_{\alpha\alpha}$,
$ u^{(\pm)}_{\beta\alpha}$
and
$\tilde{\Psi}^{(\pm)}_{\alpha}$
taking in mind for example the case of sign ``$+$''. We have from
(\ref{Aeq0}):
$$
{\cal A}=\sum\limits_{j:\, U_j\in{\cal U}_{\alpha}^d}
\psi_{\alpha}^{(j)} \langle\,\cdot\, ,\psi_{\alpha}^{(j)}   \rangle +
\Psi_{\alpha} \Psi_{\alpha}^{*} -
\sum\limits_{j:\, U_j\in{\cal U}_{\alpha}^d}
\psi_{\alpha}^{(j)} \langle\,\cdot\, , 
Q_{\alpha\beta} u_{\beta}^{(j)}  \rangle
 -\Psi_{\alpha}  u_{\beta\alpha}^{*} Q_{\alpha\beta}^{*}.
$$
It follows from the completeness relations $U^{(\pm)*}U^{(\pm)}=I-P$
for wave operators $U^{(\pm)}$ that
\be
\label{Aeq1}
 \Psi_{\alpha} u_{\beta\alpha}^{*}\equiv  
u_{\alpha\alpha} u_{\beta\alpha}^{*}=
-u_{\alpha\beta} u_{\beta\beta}^{*} -
\sum\limits_{z_j\in\sigma({{\bf H}})}
 u_{\alpha}^{(j)} \langle \,\cdot\, ,u_{\beta}^{(j)}\rangle.
\ee
Since $ u_{\beta\beta}^{*} Q_{\alpha\beta}^{*} =
(Q_{\alpha\beta} u_{\beta\beta})^{*} = u_{\alpha\beta}^{*},$ 
we can write
with a help of~(\ref{Aeq1}) that
$$
 {\cal A}=u_{\alpha\alpha} u_{\alpha\alpha}^{*} + 
u_{\alpha\beta} u_{\alpha\beta}^{*} +
    \sum\limits_{j:\, U_j\in{\cal U}_{\alpha}^d}
    \psi_{\alpha}^{(j)} \langle \,\cdot\, ,
\psi_{\alpha}^{(j)}   \rangle +
\sum\limits_{ \mbox{
\scriptsize
$\begin{array}{c}
   z_j\in\sigma_d (\mbox{{\footnotesize${\bf H}$}})\\
           U_j \not\in {\cal U}_{\alpha}^d
\end{array}$
                                         }
    }
 u_{\alpha}^{(j)}  \langle \,\cdot\, , 
Q_{\alpha\beta} u_{\beta}^{(j)}   \rangle .
$$
In the last sum, the conditions
$ z_j\in\sigma_d ({{\bf H}})$  and   
$  U_j \not\in {\cal U}_{\alpha}^d  $
mean really that we deal with any $j$ such that
$U_j\in{\cal U}_{\beta}^d$.   This follows from the equalities
$ {\cal U}_1^d \bigcup {\cal U}_2^d = {\cal U}^d $  and
$ {\cal U}_1^d \bigcap {\cal U}_2^d =  \emptyset $  (see
Theorem~\ref{ThDiscSpectrum} and Corollary~\ref{CDiscSpectrum}).
For $ U_j\in{\cal U}_{\beta}^d$, the vector 
$u_{\beta}^{(j)}$ is eigenfunction of $H_{\beta}$,
$u_{\beta}^{(j)}=\psi_{\beta}^{(j)}$, and 
$Q_{\alpha\beta} u_{\beta}^{(j)} =$ 
$Q_{\alpha\beta} \psi_{\beta}^{(j)} = u_{\alpha}^{(j)} .$
Thus, ${\cal A}$ turns in
$$
{\cal A}=u_{\alpha\alpha} u_{\alpha\alpha}^{*} + 
u_{\alpha\beta} u_{\beta\alpha}^{*}
+ \sum\limits_{z_j \in \sigma_d ({{\bf H}})}
     u_{\alpha}^{(j)}  \langle \,\cdot\, ,  
u_{\alpha}^{(j)}   \rangle
     =(U^{(\pm)}U^{(\pm)*} + P)_{\alpha\alpha}.
$$
Since $U^{(\pm)}U^{(\pm)*}+P=I$ we find
 ${\cal A}=I_{\alpha}$ and this completes 
the proof of Theorem~\ref{ThExpansion}.
\medskip

Theorem~\ref{ThExpansion} means in particular that  {\it
 part $H^{c}_{\alpha}  $ of
operator $H_{\alpha}  $   acting   in   the   invariant   subspace
corresponding to it's continuous 
spectrum $\sigma_{\alpha}^c$,  is  similar
to the operator
$
A^{(0)}_{\alpha},
$
$
H^{c}_{\alpha}  =
\Psi ^{(\pm )}_{\alpha}  A^{(0)}_{\alpha}
\tilde{\Psi }^{(\pm )*}_{\alpha}
$,
and spectrum $\sigma_{\alpha}^c$ is absolutely continuous.}
\section{INNER PRODUCT MAKING NEW HAMILTONIANS SELF-ADJOINT}
\label{InnerProduct}
We introduce now a new inner product $[\, .\, ,\, .\,]_{\alpha}$  in
$
             {{\cal H}} _{\alpha}  ,
$
$
            [f_{\alpha},g_{\alpha}]_{\alpha}= 
\langle  X_{\alpha}  f_{\alpha},g_{\alpha}   \rangle,
$
$f_{\alpha},g_{\alpha}\in{\cal H}_{\alpha},$
with $X_{\alpha}$ defined as in Theorem~\ref{ThInvariant},
$X_{\alpha}=I_{\alpha} +Q_{\alpha\beta} Q_{\alpha\beta}^{*}$,
${\alpha}=1,2.$
  The operator $X_{\alpha}  $ is  positive  definite, 
$X_{\alpha}\geq I_{\alpha}$.
This means that $[\, .\, ,\, .\,]_{\alpha}$ 
satisfies all the  axioms  of
inner product.

\begin{theorem}\label{ThInnerProduct}
The operator $H_{\alpha},  $ ${\alpha}=1,2,$  
is self-adjoint on ${\cal D}(A_{\alpha})$
in respect with the inner product $[\, .\, ,\, .\,]_{\alpha}$.
\end{theorem}

\noindent{\sl Proof.} It follows from Theorem~\ref{ThInvariant}
that operator ${{\bf H}}'$ is
self-adjoint in ${\cal H}={\cal H}_1 \oplus{\cal H}_2 $ 
in respect with the inner product
$[\, .\, ,\, .\,]$, $[f,g]=[Xf,g]$ 
with $X=\mathop{\rm diag}\{ X_1, X_2 \}$. Did, since
${\cal Q}^{-1}={\cal Q}^{*}X^{-1}=X^{-1}{\cal Q}^{*}$, 
we have for
$f,g\in{\cal D}({{\bf H}}')={\cal D}({{\bf H}})=
{\cal D}(A_1)\oplus{\cal D}(A_2)$:
$$
    [{{\bf H}}'f,g]=\langle    
X{\cal Q}^{-1} {{\bf H}} {\cal Q} f,g  \rangle =
 \langle  X\cdot X^{-1} {\cal Q}^{*}{{\bf H}}{\cal Q} f, g \rangle   =
$$
$$
   = \langle f, {\cal Q}^{*} {{\bf H}}{\cal Q} g  \rangle  =
     \langle f, X\cdot X^{-1} {\cal Q}^{*} {{\bf H}} {\cal Q} g \rangle   =
     [f,{{\bf H}}'g].
$$
Here, we used the fact that in the case of~(\ref{resolvest}),
${\cal Q} f \in {\cal D}(A_1)\oplus{\cal D}(A_2)$
if
$ f \in {\cal D}(A_1)\oplus{\cal D}(A_2).$

Taking elements $f,$ $g$ in the equality $[{{\bf H}}'f,g]=[f,{{\bf H}}'g]$
in the form $f=\{  f_1, 0  \}$, $g=\{ g_1, 0  \}$
or $f=\{ 0, f_2  \}$, $g=\{  0, g_2 \}$
with one of the components equal to zero 
and $f_{\alpha},g_{\alpha}\in{\cal D}(A_{\alpha})$,
${\alpha}=1,2$, one comes to the statement of Theorem.
\medskip

\noindent{\bf Remark.} This Theorem may be proved also in another way
making use of the equality
\be
\label{QQ}
               I_{\alpha} + Q_{\alpha\beta} Q_{\alpha\beta}^{*}  =
           \sum\limits_{j:\,\, U_j\in{\cal U} _{\alpha}^d}
 \tilde{\psi}_{\alpha}^{(j)} \langle\,\cdot\, ,
\tilde{\psi}_{\alpha}^{(j)} \rangle +
\tilde{\Psi}^{(\pm)}_{\alpha}  \tilde{\Psi}^{(\pm)} {}^{*}_{\alpha}
\ee
which is valid for both signs ``$+$'' and ``$-$''.
In this case, a self-adjointness of $H_{\alpha}$ in respect with
$[\,\cdot\, , \,\cdot\, ]_{\alpha}$  follows
from the fact that it's spectrum is real and also from
relations
$
      H^{*}_{\alpha}  \tilde{\Psi }^{(\pm )}_{\alpha}
    \tilde{\Psi }^{(\pm )*}_{\alpha}  =
$
$
      \tilde{\Psi }^{(\pm )}_{\alpha}  A^{(0)}_{\alpha}
          \tilde{\Psi }^{(\pm )*}_{\alpha}  =
$
$
               \tilde{\Psi }^{(\pm )}_{\alpha}
         \tilde{\Psi }^{(\pm )*}_{\alpha}  H_{\alpha}
$.
The equality~(\ref{QQ}) itself is proved by calculating it's right part
in the same way as it was done where the completeness relations
(\ref{completeness}) were established
(see proof of Theorem~\ref{ThExpansion}).
\section{SCATTERING PROBLEM}\label{ScatteringProblem}
We establish now that operators
$
\Psi_{\alpha}^{(+)}
$
and
$
\Psi_{\alpha}^{(-)}
$
play the same important role describing a time asymptotics of solutions
of the Schr\"{o}dinger equation
\be
\label{Schrod}
i\frac{d}{dt}f_{\alpha}(t)=H_{\alpha} f_{\alpha}(t)
\ee
 as in the usual self-adjoint case~\cite{LadyzhFaddeev},
\cite{Faddeev64}.

\begin{theorem}\label{ThWave}
Operator $U_{\alpha}(t)=\exp(iH_{\alpha} t)\exp(-iA_{\alpha}^{(0)} t)$
converges strongly if $t\rightarrow\mp\infty$,
in respect with the norm
$\parallel\cdot\parallel_{\alpha}^X$
corresponding to the inner product 
$[\,\cdot\, ,\,\cdot \,]_{\alpha}$
in ${{\cal H}} _{\alpha}$. The limit is equal to
$\,\,\, s\!-\!\mathop{\rm lim}\limits_{t\rightarrow\mp\infty} 
U_{\alpha}(t)=\Psi_{\alpha}^{(\pm)}.$
\end{theorem}

Since the norms $\| \,\cdot\, \|^{X}_{\alpha} $ and 
$\| \,\cdot\,\| $ in ${\cal H}_{\alpha}$
are equivalent,
$\|f\|\leq\|f\|^{X}_{\alpha} \leq 
(1+\|Q_{\alpha\beta}\|\cdot\|Q_{\beta\alpha}\|)^{1/2}\|f\|$,
the same statement takes place also in respect 
with the initial norm
$\| \,\cdot\,\| $.

\begin{theorem}\label{ThScattering}
For any element $f_{\alpha}^{(-)}\in{\cal H}_{\alpha}^c$
one can find such unique element  $f_{\alpha}^{(0)}$
that solution
$f_{\alpha}(t)=\exp(-iH_{\alpha} t)f_{\alpha}^{(0)}$
of Eq.~(\ref{Schrod})
satisfies the asymptotic condition
$$
\mathop{\rm lim}\limits_{t\rightarrow-\infty}
\parallel{f_{\alpha}(t)-\exp(-iA_{\alpha}^{(0)}t)
f_{\alpha}^{(-)}}\parallel_{\alpha}^X=0.
$$
There exists the unique element 
$f_{\alpha}^{(+)}\in{\cal H}_{\alpha}^c$ such that
$$
\mathop{\rm lim}
\limits_{t\rightarrow+\infty}\parallel{f_{\alpha}(t)-
\exp(-iA_{\alpha}^{(0)}t)
f_{\alpha}^{(+)}}\parallel_{\alpha}^X=0.
$$
Elements $f_{\alpha}^{(-)}$ and $f_{\alpha}^{(+)}$ are 
connected by the relation
$f_{\alpha}^{(+)}=S^{({\alpha})}f_{\alpha}^{(-)}$
with
\newline
$
     S^{({\alpha})}=
$
$
      \Psi_{\alpha}^{(-)-1}\Psi_{\alpha}^{(+)}=
$
$
       \tilde{\Psi}_{\alpha}^{(-)*}\Psi_{\alpha}^{(+)}=
$
$
         \Psi_{\alpha}^{(-)*} X_{\alpha} \Psi_{\alpha}^{(+)}.
$
\end{theorem}

We do not give here proofs of the Theorems~\ref{ThWave}
and~\ref{ThScattering} because they are exactly
the same as in the case of one-particle Schr\"{o}dinger operator in
Ref.~\onlinecite{Faddeev63}

Theorem~\ref{ThScattering}  gives the 
non-stationary formulation of the scattering problem
for a system described by Hamiltonian $H_{\alpha}$.
Moreover $S^{({\alpha})}$ is a scattering   operator for this system.

\begin{theorem}\label{ThScChannel}
Scattering operator $S^{({\alpha})}$ coincides with the component
$s_{\alpha\alpha}$ of the scattering operator $S$, $S=U^{(-)*}U^{(+)}$,
for a system described by the two--channel Hamiltonian ${{\bf H}}$.
\end{theorem}

\noindent{\sl Proof.} Let us show that operator $S^{({\alpha})}$
has the kernel $s_{\alpha\alpha}(\lambda,\lambda')$ 
given by Eq.~(\ref{scattering}).
To do this, remember that 
$\tilde{\Psi}_{\alpha}^{(-)}=\Psi_{\alpha}^{(-)}-
Q_{\alpha\beta} u_{\beta\alpha}^{(-)}$
(see Theorem~\ref{ThExpansion}). Therefore,
$$
S^{({\alpha})}=(\Psi_{\alpha}^{(-)*}-
u_{\beta\alpha}^{(-)*} Q_{\alpha\beta}^{*})\Psi_{\alpha}^{(+)}=
\Psi_{\alpha}^{(-)*} \Psi_{\alpha}^{(+)} + 
u_{\beta\alpha}^{(-)*} Q_{\beta\alpha} \Psi_{\alpha}^{(+)} =
u_{\alpha\alpha}^{(-)*} u_{\alpha\alpha}^{(+)} +  
u_{\beta\alpha}^{(-)*} u_{\beta\alpha}^{(+)} .
$$
Here, we have used the properties
${\Psi}^{(\pm)}_{\alpha} = u_{\alpha\alpha}^{(\pm)} $,
$Q_{\alpha\beta}^{*}=-Q_{\beta\alpha}$ and
$Q_{\beta\alpha} \Psi_{\alpha}^{(+)} =u_{\beta\alpha}^{(+)}$ 
established above.
Since
$$
      u_{\alpha\alpha}^{(-)*} u_{\alpha\alpha}^{(+)} +  
      u_{\beta\alpha}^{(-)*} u_{\beta\alpha}^{(+)} =
\left(  U^{(-)*}U^{(+)}  \right) _{\alpha\alpha} = s_{\alpha\alpha},
$$
we come to the statement of Theorem. The proof of
Theorem~\ref{ThScattering} is completed.
\medskip

A kernel of the scattering operator $S^{({\alpha})}$ may be presented
also in a usual way~(\ref{scattering}) in terms of the $t$-matrix
$t_{\alpha}(z)=W_{\alpha}-W_{\alpha}(H_{\alpha}-z)^{-1}W_{\alpha}$, 
taken on the energy--shell.
Note that $t_{\alpha}(z)$ differs from $T_{\alpha\alpha}(z)$ 
introduced in Sec.~\ref{Initial}.
Did, easy calculations show that
\be
\label{Tmatrix}
t_{\alpha}(z)=B_{\alpha\beta} 
[I_{\beta} +Q_{\beta\alpha} (A_{\alpha} -z)^{-1} B_{\alpha\beta}]^{-1} 
Q_{\beta\alpha} .
\ee
Using the basic equation~(\ref{QbasicSym}) one can rewrite~(\ref{Tmatrix})
in the form
$$
    t_{\alpha}(z)=T_{\alpha\alpha}(z) + \tilde{t}_{\alpha}(z)
$$
where
$$
\tilde{t}_{\alpha}(z)= B_{\alpha\beta} [A_{\beta} -
B_{\beta\alpha}(A_{\alpha} -z)^{-1}
B_{\alpha\beta}]^{-1}Q_{\alpha\beta} (A_{\alpha} -z)\not\equiv 0.
$$
However the additional term $\tilde{t}_{\alpha}(z)$
is evidently disappearing on the energy-shell due to presence of the
difference $A_{\alpha} -z$ as an end factor. Actually, in the
diagonal representation~(\ref{Neumann}),(\ref{multi}),
$A_{\alpha} -z$ acts as the factor $\lambda-z$ 
vanishing at $z=\lambda+i0$.
Therefore, kernels of $t$-matrices 
$t_{\alpha}$ and $T_{\alpha\alpha}$ coincide
on the energy surface.

Note also that in our case
$\sigma_1^c\bigcap\sigma_2^c=\emptyset$. Hence we have
$s_{{\beta}{\alpha}}=0$ and 
$S^{({\alpha})}=s_{{\alpha}{\alpha}}$ is unitary.
\medskip
	
\acknowledgements
Author is thankful to Prof. B.S.Pavlov and Dr.~K.A.Makarov for
support and interest to this work.  The author is much indebted to
Prof.~R.A.Minlos and  participants of his seminar in Moscow
University for discussion of the presented results.  Particularly,
the author is grateful to Dr.~S.A.Stepin for the literary
indications.

{The work is supported in part by the International Science Foundation
(Grants~\#RFB000 and~\#RFB300) and the RAS Academy of Natural Sciences.}

\newpage

\end{document}